\documentclass[12pt]{amsart}
\usepackage[top=30mm,bottom=30mm,left=30mm,right=30mm]{geometry}

\usepackage{color}
\usepackage{lscape}
\usepackage{xypic}

\newcommand{\la}{\lambda}

\newcommand{\C}{\mathbb{C}}

\newcommand{\Z}{\mathbb{Z}}

\newcommand{\OG}{\mathrm{OG}}
\newcommand{\LG}{\mathrm{LG}}
\newcommand{\SO}{\mathrm{SO}}
\newcommand{\QH}{{QH}}
\newcommand{\Sp}{\mathrm{Sp}}

\newcommand{\cG}{\mathcal{G}}

\renewcommand{\mod}{\mathrm{mod}\;}

\renewcommand{\u}{\pmb{u}}
\renewcommand{\v}{\pmb{v}}
\newcommand{\e}{\pmb{e}}

\newcommand{\diag}{{\mathrm{diag}}}

\theoremstyle{definition}
\newtheorem{prop}{Proposition}[section]
\newtheorem{thm}{Theorem}[section]
\newtheorem{lem}[prop]{Lemma}%[section]
\newtheorem{cor}[prop]{Corollary}%[section]

\theoremstyle{remark}
\renewcommand{\theequation}{\arabic{section}.\arabic{equation}}

\begin{document}
\title[Factorial $P$- and $Q$-Schur functions]
{Factorial $P$- and $Q$-Schur functions represent equivariant quantum Schubert classes}
% and equivariant quantum cohomology of 
%maximal isotropic Grassmannians}

\begin{abstract}
We find presentations by generators and relations for the equivariant 
quantum cohomology rings of the maximal isotropic Grassmannians of types B,C and D, and we find polynomial representatives for the Schubert classes in these rings. These representatives are given in terms of the same Pfaffian formulas which appear in the theory of factorial $P$- and $Q$-Schur functions. After specializing to equivariant cohomology,  we interpret the resulting presentations and Pfaffian formulas in terms of Chern classes of tautological bundles. 
%To discuss the type D Orthogonal Grassmannians, we investigate the algebraic properties of the Dynkin involution, and provide a (quantum) cohomological interpretation of the odd limit of factorial $P$-Schur functions, which might be of independent interest.
\end{abstract}

\author{Takeshi Ikeda} \address{Department of Applied Mathematics,                                            
Okayama University of~Science, Okayama 700-0005, Japan} \email{ike@xmath.ous.ac.jp}

\author{Leonardo~C.~Mihalcea}
\address{Department of Mathematics, Virginia Tech, 460 McBryde, Blacksburg, VA 24061, USA}
\email{lmihalce@math.vt.edu}

\author{Hiroshi Naruse}\address{Graduate School of Education, Okayama University, Okayama 700-8530, Japan}
\email{rdcv1654@cc.okayama-u.ac.jp}

%\author{Takeshi Ikeda}
%\author{Leonardo C.~ Mihalcea}
%\author{Hiroshi Naruse}
%\date{2011 Dec. 05}
\date{February 4, 2014}


%\thanks{}

\maketitle

%\tableofcontents

\section{Introduction} 
In this paper, we study the torus equivariant quantum cohomology rings of the maximal isotropic Grassmannians in classical types B,C and D.
The main purpose is to find a presentation for the ring and prove a Giambelli formula for the equivariant quantum 
Schubert classes. 

For $n$ fixed, we denote the {\em Lagrangian Grassmannian\/} by $\LG(n)$ in type C, 
which parametrizes subspaces in a symplectic vector space $\C^{2n}$ which are {\em Lagrangian\/},
i.e. the subspaces of dimension $n$ which are isotropic with respect to a symplectic form.
By $\OG(n)$, we denote the maximal {\em Orthogonal Grassmannian\/} 
in type D, which parametrize $(n+1)$-dimensional isotropic subspaces of an orthogonal vector space $\C^{2n+2}$.
~(We only consider one of the two connected components in type D, see \S \ref{ss:D} below.)
For type B, we consider the maximal isotropic Grassmannian for the odd dimensional orthogonal space $\C^{2n+1}$.
The Grassmannian of type B is  
known to be isomorphic to $\mathrm{OG}(n)$ as algebraic varieties. Although the tori acting on these varieties are different,
we can deduce results for type B from those of type D, see \S \ref{ssec:BD} for details.

Let $\cG_n$ denote one of $\LG(n)$ or $\OG(n)$ and let $T$ be the maximal torus of the complex symplectic group $\Sp_{2n}$, respectively the complex special orthogonal group $\SO_{2n+2}$. Denote by $S:=H^*_T(pt)$ the integral equivariant cohomology of a point, which is the polynomial ring $\Z[t_1,\ldots,t_n]$ (respectively $\Z[t_1, \ldots, t_{n+1}]$) in the characters of $T$. The equivariant quantum cohomology ring $\QH^*_T(\cG_n)$, defined for more general varieties by Kim \cite{kim:oneqq}, 
is a graded $S[q]$-algebra, where the quantum parameter $q$ has (complex) degree $n+1$ or $2n$ respectively.
It has an $S[q]$-basis consisting of Schubert classes $\sigma_\lambda$, where $\lambda=(\lambda_1> \cdots> \lambda_k>0)$ varies in the set $\mathcal{SP}(n)$ of strict partitions included in the staircase $(n, n-1, \ldots, 1)$. 
The multiplication \[ \sigma_\lambda * \sigma_\mu = \sum_{d \ge 0, \nu \in \mathcal{SP}(n)} c_{\lambda,\mu}^{\nu,d} q^d \sigma_\nu \] is determined by the ($3$-point, genus $0$) equivariant Gromov-Witten (GW) invariants, defined by Givental \cite{givental:eqgw}.
The coefficients $c_{\lambda,\mu}^{\nu,d}$ are homogeneous polynomials in $S$, and those $c_{\lambda,\mu}^{\nu,d}$ of polynomial degree $0$
are non-negative integers equal to the ordinary GW invariants counting rational curves of degree $d$ passing through general translates of Schubert varieties. If $d=0$, $c_{\lambda,\mu}^{\nu,d}$ is a structure constant of the equivariant cohomology ring $H^*_T(\cG_n)$.~There is a $\Z$-algebra isomorphism $\QH^*_T(\cG_n)/\langle S_+ \rangle \simeq \QH^*(\cG_n)$ to the quantum cohomology of $\cG_n$, and an $S$-algebra isomorphism $\QH^*_T(\cG_n)/\langle q \rangle \simeq H^*_T(\cG_n)$ to the equivariant cohomology ring; here $S_+$ consists of the elements in $S$ of positive degrees.

The main goal of this paper is to solve the {\em Giambelli problem} for $\QH^*_T(\cG_n)$: ~(1) we find a presentation with generators and relations of $\QH^*_T(\mathcal{G}_n)$;~(2) we identify a set of polynomials in the given generators which are sent to Schubert classes $\sigma_\lambda$ (the {\em Giambelli formula}). 
It turns out that a natural combinatorial framework for these statements is given by the {\em factorial $P$- and $Q$- Schur functions} $P_\lambda(x|t), Q_\lambda(x|t)$. These functions are a slight variation of those introduced by Ivanov \cite{Iv} - see \S \ref{ssec:deffactP} below -  and they are deformations of the ordinary $P$- and $Q$-Schur functions $P_\lambda(x), Q_\lambda(x)$ (see Schur's paper \cite{Sc}).
It is a general feature of the theory of $P$- and $Q$- Schur functions that they can be expressed as Pfaffians of skew-symmetric matrices, and this extends to their factorial deformation \cite{Iv}. The relation between Ivanov's factorial $P$- and $Q$- Schur functions and geometry of $\cG_n$ was first established by the first and third authors of this paper in \cite{Ik,EYD}, where they solve the Giambelli problem for the equivariant ring $H^*_T(\cG_n)$ using these functions. Therefore, one expects that some $q$-deformations of these Pfaffians represent Schubert classes in $\QH^*_T(\cG_n)$. The pleasant - and somewhat surprising - fact proved in this paper is that the same, undeformed Pfaffian formulas answer the Giambelli problem in the quantum ring. This is reminiscent of the situation regarding the (equivariant) quantum Giambelli problem for the type A Grassmannian \cite{bertram:quantum, Mihalcea-Giambelli}, and the quantum Giambelli problem for $\cG_n$ \cite{KT-LG, KT-OG}.

\subsection{Statement of results} To state our results precisely, we fix some notation. Let $x=(x_1,x_2,\ldots)$ be an infinite sequence of variables.
Let $P_i(x),Q_i(x)$ denote Schur's $P$- and $Q$-functions
for partitions with one part (cf. \cite{Mac} III.8).~Set $\Gamma'=\Z[P_1(x),P_2(x),\ldots]$ 
and $\Gamma=\Z[Q_1(x),Q_2(x),\ldots]$ and recall that these rings have a $\Z$-basis 
given by the $P$- and $Q$-Schur functions $P_\lambda(x)$ and $Q_\lambda(x)$, where 
$\lambda=(\lambda_1 > \cdots > \lambda_k >0)$ varies in the set of strict partitions $\mathcal{SP}$. In fact, $\Gamma$ is a subring of $\Gamma'$, because $Q_\lambda(x) = 2^{\ell(\lambda)} P_\lambda(x)$, where $\ell(\lambda) = k$ denotes the length of $\lambda$.~Let $t=(t_1,t_2,\ldots)$ be a sequence of indeterminates
and set $\Z[t]:=\Z[t_1,t_2,\ldots]$. Then $\Z[t]\otimes_\Z\Gamma'$ is naturally 
a graded ring with $\deg P_i(x)=i,\deg t_i=1$,
and $\Z[t]\otimes_\Z\Gamma$ is a graded subring of it. 
For each $\la \in \mathcal{SP}$, 
the corresponding factorial $P$-function 
$P_\la(x|t)$ (resp. $Q$-function $Q_\la(x|t)$) is a homogeneous 
element of $\Z[t]\otimes_\Z\Gamma'$ (resp. $\Z[t]\otimes_\Z\Gamma$).
If we set all the parameters $t_i$ to zero, then
$P_\la(x|t)$ specializes to the ordinary 
$P$-Schur function $P_\la(x),$
and similarly for $Q_\la(x|t)$. 
We use the convention that $t_i=0$ if $i>n$ for $\LG(n)$, respectively if $i > n+1$ for $\OG(n)$.
The following is the main result of this paper.

\begin{thm}\label{thm:main} (a)  There is an isomorphism of 
graded $S[q]$-algebras \[
S[q][P_1(x|t),\ldots,P_{2n}(x|t)]/I_n^{(q)}
\longrightarrow \QH_T^*(\OG(n)),
\]
where 
$I_n^{(q)}$ is the ideal generated by 
$P_{n+1}(x|t),\ldots,P_{2n-1}(x|t), P_{2n}(x|t)+(-1)^n q$.
Moreover, the image of $P_\la(x|t)\;(\la \in \mathcal{SP}(n))$
is the Schubert class $\sigma_\la.$ 

(b) There is an isomorphism of 
graded $S[q]$-algebras
\begin{equation*}
S[q][Q_1(x|t), \ldots, Q_n(x|t), 2Q_{n+1}(x|t), Q_{n+2}(x|t), \ldots, Q_{2n}(x|t)]/J_n^{(q)} \longrightarrow \QH_T^*(\LG(n)) \/, 
\end{equation*}

\noindent where 
$J_n^{(q)}$ is the ideal generated by 
$2Q_{n+1}(x|t)-q,Q_{n+2}(x|t), \ldots,Q_{2n}(x|t)$.
Moreover, the image of $Q_\la(x|t)\;(\la \in \mathcal{SP}(n))$
is the Schubert class $\sigma_\la$. 
\end{thm}

In Theorems \ref{thm:OG} and \ref{thm:LG} below we also give a presentation with (independent) generators and relations. In that case, the ideal of relations contains the quadratic identities which are satisfied by the (factorial) $P$- and $Q$-functions (cf. (\ref{eq:quadratic}), (\ref{E:Qqvanishing})).%

The specialization at $q=0$ in the relations from (a) and (b) recover the equivariant cohomology rings. The Pfaffian formulas and the (specialized) quadratic identities have a geometric interpretation in terms of equivariant Chern classes of the tautological bundles on $\OG(n)$, explained in \S \ref{ss:chernOG} below.  This extends authors' work in \cite{DSP}, where it was done for the Lagrangian Grassmannian. The Pfaffian formula 
was first proved by Kazarian \cite{kazarian}
in the context of degeneracy loci formulas of vector bundles.
Our proof is different, and it is based on our earlier results from \cite{DSP}, in which we introduced the double Schubert polynomials, a canonical family of
polynomials identified with the Schubert classes in the equivariant cohomology of full flag manifolds of types B,C,D. 
Kazarian's method was recently employed by Anderson and Fulton \cite{anderson.fulton} to extend the single Pfaffian expression of the double Schubert polynomials to a wider class of Weyl group elements called {\em vexillary signed permutations\/}. It should be noted that Tamvakis \cite{tamvakisG/P, tamvakis-survey} proved a general combinatorial formula that expresses the double Schubert polynomials as explicit positive linear combinations of products of Jacobi-Trudi determinants times (at most) a single Schur Pfaffian. 

\subsection{Related work and idea of proof} The connection between the (ordinary) $P$- and $Q$-Schur functions and the Giambelli problem for the ordinary cohomology of $\OG(n)$ and $\LG(n)$ was discovered by Pragacz \cite{pragacz}; this is analogous to the classical fact that the ordinary Schur functions solves the Giambelli problem for the type A Grassmannian - see e.g. \cite{fulton:young}. Kresch and Tamvakis \cite{KT-LG,KT-OG} used intersection theory on the Quot  schemes in classical types - an argument simplified later by Buch, Kresch and Tamvakis \cite{BKT:gw} - to find a presentation for the quantum cohomology rings of $\LG(n)$ and $\OG(n)$. They proved that the Pfaffian formula for the $P$- and $Q$-Schur functions gives representatives for the quantum Schubert classes. The Quot scheme approach was pioneered by Bertram \cite{bertram:quantum}, who proved that the determinantal formula for the Schur functions gives the quantum Giambelli formula in Witten's presentation
 \cite{witten} of the quantum cohomology ring of the type A Grassmannian. There is large body of literature dedicated to the (non-equivariant) quantum Giambelli problem for other homogeneous spaces - see e.g. \cite{AS,BKT:qiso, ciocan:partial,FGP,kim:toda} and references therein. The equivariant quantum Giambelli formula for partial flag manifolds was recently and independently solved by Anderson and Chen \cite{anderson.chen:qdouble} (by using Quot schemes), and by Lam and Shimozono \cite{lam.shimozono:qdsp} (using the ``Chevalley approach" explained in the next paragraph). The answer was given in terms of specializations of Fulton's universal Schubert polynomials \cite{fulton:universal}; for the full flag manifolds this specialization recovers the quantum double Schubert polynomials which appeared in a paper by Kirillov and Maeno \cite{kirillov.maeno:qdouble}.

The proof of theorem \ref{thm:main} is logically independent on earlier results from \cite{BKT:gw, KT-LG,KT-OG} regarding quantum cohomology of $\LG(n)$ and $\OG(n)$, and in fact our methods give an alternate proof of those results. We rely on the characterization of the equivariant quantum cohomology ring of any homogeneous variety in terms of the {\em Chevalley formula} - see Theorem \ref{prop:chareqq} below. This was proved by the second author 
in \cite{Mihalcea-EQChevalley} (initially in \cite{Mihalcea-EQSC} for Grassmannians) and it was succesfully used to solve the Giambelli problem for the equivariant quantum cohomology ring of the Grassmannian \cite{Mihalcea-Giambelli} and partial flag manifolds \cite{lam.shimozono:qdsp}.~We show that the product $P_{1}(x|t) \cdot P_\lambda(x|t) $ satisfies the Chevalley formula in the equivariant quantum ring, modulo the given ideal.~The proof uses a Gr\"obner basis argument showing that the images of factorial $P$-functions form an $S[q]$-basis for the given quotient ring.
By the aforementioned characterization theorem, this gives the result. A similar approach works for $\LG(n)$, although there are some technical differences.
%due to some multiples of $2$ which appear in the type C. 

\subsection{Organization} In Section \ref{FacP}, we present preliminary results on the factorial $Q$- and $P$-functions.
In Section \ref{sec:ecoh}, we fix some notations for the maximal isotropic Grassmannians
and state the characterization results for the equivariant quantum cohomology rings by the Chevalley rule.
In sections \ref{sec:OG} and \ref{sec:LG} we prove the main theorem respectively for 
$\OG(n)$ and $\LG(n)$. In 
Appendix \ref{ss:chernOG}, we discuss the (non-quantum) equivariant cohomology of $\OG(n)$, and we give a geometric interpretation, in terms of Chern classes, of the algebraic quantities from our main theorem.
%interpret various al
%We will show how the algebraic objects 
%entered in the previous arguments are written in terms of more geometric languages.

\setcounter{equation}{1}
\section{Factorial $P$-Schur and $Q$-Schur functions}\label{FacP}

The goal of this section is to recall the definition and the main properties of factorial $P$-Schur and $Q$-Schur functions, following mainly Ivanov's paper \cite{Iv}.

\subsection{Definition of factorial $P$-functions}\label{ssec:deffactP}
Let 
$\la=(\la_1>\cdots>\la_k>0)$
be a strict partition with $k \le N$; the quantity $\ell(\lambda):=k$ is 
{\em the length} of $\la$. Let $t= (t_1, t_2, \ldots)$ be a sequence of indeterminates.
The {\it generalized 
factorial\/} is defined by 
$$
(x|t)^k=(x-t_1)\cdots(x-t_k).
$$
Following \cite{Iv}, define the {\it factorial $P$-Schur function\/}
$P_\la(x_1,\ldots,x_N|t)$ to be
$$
P_\la(x_1,\ldots,x_N|t)=\frac{1}{(N-\ell(\la))!}
\sum_{w\in S_N}w
\left(
\prod_{i=1}^{\ell(\la)}(x_i|t)^{\la_i}
\prod_{i=1}^{\ell(\la)}\prod_{j=i+1}^N
\frac{x_i+x_j}{x_i-x_j}
\right),
$$
where $w\in S_N$ permutes the variables $x_1,\ldots,x_N$. The polynomials $P_\la(x_1,\ldots,x_N|t)$ are not stable when the number of variables increases, i.e.,
$P_\la(x_1,\ldots,x_N,0,|t) \neq P_\la(x_1,\ldots,x_N|t)$ in general. However, these polynomials satisfy the weaker stability property $$
P_\la(x_1,\ldots,x_N,0,0|t) = P_\la(x_1,\ldots,x_N|t).$$
Therefore it makes sense to define two projective limits \[ P_\lambda^+(x|t):= \varprojlim P_\la(x_1,\ldots,x_{2N}|t); % removed "still" 
\quad P_\lambda^-(x|t) =\varprojlim P_\la(x_1,\ldots,x_{2N+1}
|t) \]  taken over even, respectively odd, number of variables. 

\bigskip

\noindent{Remark 1.} It is known that $P_\lambda^{+}(x|t)$ corresponds to Schubert classes equivariant under a torus $T$ (see e.g.~\S\ref{ss:chernOG} below).
The odd limit $P_\lambda^-(x|t)$ also has a geometric relevance. We will study in more detail the function $P_\lambda^-(x|t)$
in an upcoming paper.

\bigskip

In what follows we only consider the even limit, denote it by $P_\lambda(x|t):= P_\lambda^+(x|t)$, and we refer to it as {\em the factorial P-Schur function}.
If $t_i =0$ for $i \ge 1$ then one recovers the definition of the ordinary $P$-Schur function defined in \cite[Ch. III, \S 8]{Mac}. Recall that $\Gamma'$ denotes the ring $\Z[P_1(x), P_2(x), \ldots]$. Then $P_\lambda(x|t)$ is an element of
$\Z[t]\otimes_\Z\Gamma'$.
The {\it factorial $Q$-Schur function\/} is defined by 
\begin{equation}
Q_\la(x|t)=2^{\ell(\la)}P_\la(x|0,t_1,t_2,\ldots).
\label{def:Q}
\end{equation} Again, if all $t_i$'s equal to $0$ then one recovers the ordinary $Q$-Schur function, defined by Schur \cite{Sc} in relation to projective representations of the symmetric group; see \cite[Ch. III, \S 8]{Mac} for more on $Q$-Schur functions. Set $\Gamma:= \Z[Q_1(x), Q_2(x), \ldots]$. Then the element $Q_\lambda(x|t) \in  \Z[t]\otimes_\Z \Gamma$ is called the factorial $Q$-Schur function.
By \cite[Ch. III, \S 8, (8.9)]{Mac} the ordinary $P$-Schur, resp. $Q$-Schur functions form a $\Z$-basis of $\Gamma'$, resp. $\Gamma$. This implies that the functions
$P_\la(x|t)$ and $Q_\la(x|t)$, when $\lambda$ varies over strict partitions, form bases over $\Z[t]$ for the rings $\Z[t]\otimes_\Z \Gamma'$ respectively $\Z[t]\otimes_\Z \Gamma$.

\subsection{Quadratic identities, the Pfaffian formula, and recurrences}

It is easy to see from the definition that $P_1(x|t) = P_1(x)=\sum_{i=1}^\infty x_i$ and more generally
 \begin{equation}\label{E:recPi} P_i(x|t) = \sum_{j=0}^{i-1} (-1)^j e_j(t_1, \ldots, t_i) P_{i-j}(x), \/ \end{equation} where $e_j(t_1, \ldots, t_i)$ is the $j$th elementary symmetric function. The functions $P_{k,l}(x|t)$ are determined by $P_i(x|t)$ using certain quadratic identities, which will play a key role in this paper. To define these relations, let $x=(x_1,\ldots,x_a)$ and $y
=(y_1,\ldots,y_b)$ be two sets of indeterminates, and set 
$$h_k(x|y)=\sum_{i+j=k}h_i(x_1,\ldots,x_a)
e_j(y_1,\ldots,y_b)$$ where $h_i(x_1,\ldots,x_a)$ denotes the $i$th complete homogeneous symmetric functions.
%; see \cite{St} for
%more information on supersymmetric functions.
 
\begin{prop}\label{prop:quadratic} Let $k,l$ be positive integers such that $k\geq l.$ Then 
$$
P_{k,l}(x|t)=P_k(x|t)P_l(x|t)+
\sum_{(r,s)}g_{k,l}^{r,s}(t)P_r(x|t)P_s(x|t),
$$
where the sum is over the set $\mathcal{I}_{k,l}$ consisting of pairs $(r,s)$ such that
$(r,s)\neq (k,l)$ and   
\begin{equation}
k\leq r\leq k+l,\;
0\leq s\leq l,\;
r+s\leq k+l,\label{eq:rs}
\end{equation} 
and where
\begin{equation}\label{E:g}
g_{k,l}^{r,s}(t)=\begin{cases}
(-1)^{l-s}2h_{k+l-r-s}(t_{k+1},\ldots,t_{r+1}|t_{s+2},\ldots,t_l) & \mbox{if}\;s\geq 1\\
(-1)^{l-s}h_{k+l-r-s}(t_{k+1},\ldots,t_{r+1}|t_{1},\ldots,t_l) & \mbox{if}\;s=0
\end{cases}.
\end{equation}
\end{prop}
\noindent{Proof.} For the factorial $Q$-functions (i.e. when $t_1=0$), 
this formula was proved in
\cite{Ik}, Proposition 7.1. 
The proof from \cite{Ik} uses
equations (8.2) and (8.3) in \cite{Iv}, which have a straightforward generalization to the case of $t_1 \neq 0$.
Then the proof from {\em loc.~cit.} extends to this case. \qed

\medskip

%The quadratic relations have a geometric interpretation given in terms of identities among the Chern classes of the tautological bundles %on $\OG(n)$; in fact they can be deduced from this; see \S \ref{ss:chernOG} below. 
Since $P_{i,i}(x|t) = 0$ by \cite[Proposition 2.6 (c)]{Iv} we obtain:

\begin{cor}\label{cor:quadraticvanishing} The following identity holds for each $i \ge 1$: \begin{equation}
P_{i}(x|t)^2+\sum_{(r,s) \in \mathcal{I}_{i,i}}g_{i,i}^{r,s}(t)
P_r(x|t)P_s(x|t)=0 \/.\label{eq:quadratic}
\end{equation}
\end{cor}

The quadratic relations have a geometric interpretation given in terms of identities among the Chern classes of the tautological bundles on $\OG(n)$; in fact they can be deduced from this - see \S \ref{ss:chernOG} below. Similar quadratic recurrences and identities hold for the factorial $Q$-Schur functions, and were proved by the first author in \cite[Proposition 7.1]{Ik}.

\begin{prop}[\cite{Ik}, Proposition 7.1]
\begin{equation}\label{E:Qrec}
Q_{k,l}(x|t)=Q_k(x|t)Q_l(x|t)+
\sum_{(r,s) \in \mathcal{I}_{k,l} }f_{k,l}^{r,s}(t)Q_r(x|t)Q_s(x|t) \end{equation}
where $
f_{k,l}^{r,s}(t)
=(-1)^{l-s}2h_{k+l-r-s}(t_k,\ldots,t_r|t_{s+1},\ldots,t_{l-1})
$.
In particular, we have the identity  \begin{equation}\label{E:Qqvanishing}
Q_{i}(x|t)^2+\sum_{(r,s) \in \mathcal{I}_{i,i}}f_{i,i}^{r,s}(t)
Q_r(x|t)Q_s(x|t)=0 \/. 
\end{equation}
\end{prop}
The factorial $P$-Schur functions $P_\lambda(x|t)$ for an arbitrary strict partition $\lambda$ can be calculated recursively starting from partitions having at most two parts. To state this precisely, let $\lambda$ be a strict partition, and set 
\begin{equation} \label{E:r} r=2\left[(\ell(\la)+1)/2\right] \end{equation}
and $\la_{r}=0$ if $\ell(\la)$ is odd. 
%{\color{blue} 
In this case we make $P_{\lambda_i, 0}(x|t) = P_{\lambda_i}(x|t)$ by convention. Also by convention $P_{k,l}(x|t) = - P_{l,k}(x|t)$. This is consistent with the fact that one can define factorial $P$-Schur functions for {\em any} partition $\lambda$ (strict or not), but if two parts of $\lambda$ are equal then the corresponding function vanishes \cite[Proposition 2.6]{Iv}.
%}

\begin{prop}[\cite{Iv}] We have 
\begin{equation}
P_\la(x|t)=\mathrm{Pf}(P_{\la_i,\la_j}(x|t))_{1\leq i<j\leq r}
\label{eq:Pfaff-P}
\end{equation}
where $\mathrm{Pf}(A)$ denotes the Pfaffian of the (skew-symmetric) matrix $A$. In particular, we have the following 
recurrence relations:
\begin{equation}
P_\la(x|t)=\sum_{i=2}^{r}(-1)^iP_{\la_1,\la_i}(x|t)
P_{\la_2,\ldots,\widehat{\la_{i}},\ldots,\la_{r}}(x|t).
\label{eq:rec}\end{equation}
\end{prop}

Note that if $\ell(\la)$ is odd then (\ref{eq:rec}) is equivalent to the following identity: 
\begin{equation}
P_\la(x|t)=\sum_{i=1}^{\ell(\la)}
(-1)^{i-1}P_{\la_i}(x|t)
P_{\la_1,\ldots,\widehat{\la_{i}},\ldots,\la_{r}}(x|t).
\label{eq:rec-odd}
\end{equation}

%\bigskip
\medskip
\noindent{Remark 2.}
It is obvious from (\ref{def:Q}),
that the same formula holds for $Q_\la(x|t).$

%\medskip

\subsection{}\label{ss:ChernPk} Next we define and study some functions which have a natural interpretation in terms of Chern classes of tautological bundles. Their properties will be used in section \S \ref{ss:chernOG} below. Define $P_i^{(k)}(x|t) \in \Z[t] \otimes_Z \Gamma'$ by the generating series \begin{equation}
 1 + \sum_{i=1}^\infty 2 P_i^{(k)}(x|t)u^k = \prod_{i=1}^\infty \frac{1+x_iu}{1-x_iu} \prod_{j=1}^k (1-t_ju) \/. \label{eq:2P_i}
 \end{equation} Equation (\ref{E:recPi}) above shows that 
\begin{equation} P_k(x|t) =  P_k^{(k)}(x|t) - \frac{(-1)^k}{2} e_k(t_1, \ldots, t_k) \/. \label{eq:P=P+e/2}
\end{equation} 

\begin{lem}\label{lem:Pshit} We have
\begin{eqnarray*}P_{k+j}^{(k)}(x|t)&=&P_{k+j}(x|t)
+h_1(t_{k+1},\ldots,t_{k+j})P_{k+j-1}(x|t)
+\cdots
+h_{j-1}(t_{k+1},t_{k+2})P_{k+1}(x|t)\\
&+&h_j(t_{k+1})P_{k}(x|t)\quad (j\geq 1),\\
P_{k-j}^{(k)}(x|t)&=&P_{k-j}(x|t)
-e_1(t_{k-j+1},\ldots,t_{k})P_{k-j-1}(x|t)
+e_2(t_{k-j},\ldots,t_{k})P_{k-j-2}(x|t)
+\cdots
\\
&+&(-1)^{k-j-1}e_{k-j-1}(t_{3},\ldots,t_{k})P_{1}(x|t)
+(-1)^ke_k(t_1,\ldots,t_{k})/2\quad (j\geq 0).
\end{eqnarray*}
\end{lem}
\noindent{Proof.}
The equation is easily verified (cf. equations in \cite[p. 882]{DSP}).
\qed

\begin{lem}\label{lem:2row} We have
\begin{eqnarray}
P_{k,l}(x|t)&=&\left(P_k^{(k)}(x|t)-(-1)^ke_k(t_1,\ldots, t_k)/2\right)\left(P_l^{(l)}(x|t)+(-1)^le_l(t_1,\ldots, t_l)/2\right)
\nonumber\\
&+&2\sum_{j=1}^{l-1}(-1)^jP_{k+j}^{(k)}(x|t)P_{l-j}^{(l)}(x|t)+
(-1)^lP_{k+l}^{(k)}(x|t).
\end{eqnarray}
\end{lem} 
\noindent{Proof.} The expression of $P_{k,l}(x|t)$ in
Proposition \ref{cor:quadraticvanishing} can be
written in the above form by using Lemma \ref{lem:Pshit}.
\qed

\subsection{The Chevalley rule}

\begin{prop}[\cite{Iv},Theorem 6.2] 
Let $\la$ be a strict partition and  $r$ from equation (\ref{E:r}) above. 
Then
\begin{equation} P_1(x|t) P_\lambda(x|t) = \sum_{\mu\to \la}P_\mu(x|t) + (\sum_{i=1}^r t_{\la_i+1}) P_\la(x|t) \/, \label{eq:Chev}
\end{equation}
where $\mu\to\la$ means that $\mu$ is a strict partition obtained from $\lambda$ by adding one more box. In particular
%\[ Q_1(x|t) Q_\lambda(x|t) = \sum_{\mu\to \la}2^{\ell(\lambda) - \ell(\mu) +1} Q_\mu(x|t) + 2(\sum_{i=1}^r t_{\la_i+1}) Q_\la(x|t) \/. \] 
\[ Q_1(x|t) Q_\lambda(x|t) = \sum_{\mu\to \la}2^{\ell(\lambda) - \ell(\mu) +1} Q_\mu(x|t) + 2(\sum_{i=1}^{\ell(\lambda)} t_{\la_i}) Q_\la(x|t) \/.\]

\end{prop}

\setcounter{equation}{1}
\section{Equivariant quantum cohomology of maximal isotropic Grassmannians }\label{sec:ecoh}
The goal of this section is to fix notations for the maximal isotropic Grassmannians, and to recall the definition and some basic facts of their equivariant quantum cohomology ring.

\subsection{Schubert classes and equivariant cohomology}\label{ssec:eqqc} We recall next some basic facts about the torus equivariant cohomology ring; see \cite{brion:eq,brion:eqpoincare} for details. Although the definitions here make sense for any homogeneous space $X$, we will restrict to the case when $X=G/P$ is a homogeneous space where $G$ is a complex semisimple Lie group and $P$ is a {\em maximal} parabolic subgroup. (We will soon specialize further both $G$ and $P$.) Then $X$ is a smooth, complex, projective variety with an action of a maximal torus $T \subset G$ given by left multiplication. Consider the universal bundle $ET \to BT$. Then $T$ acts freely on $ET$ and one can define a (free) $T$-action on $ET\times X$ by $t\cdot (e, x) = (et^{-1}, tx)$. The (integral) equivariant cohomology of $X$, denoted $H^*_T(X)$, is the ordinary cohomology of the ``Borel mixed space" $X_T:= (ET \times X)/T$. The structure morphism $X \to pt$ gives $H^*_T(X)$ a structure of an $S$-algebra, where $S= H^*_T(pt)$. 
In fact, $S$ can be identified with the polynomial ring $\Z[t_1, \ldots, t_l]$ where $\{t_1,\ldots,t_l\}$ is a basis of 
the characters group (written additively) of the torus $T$. 
Each irreducible, closed, subvariety $Y \subset X$ of (complex) codimension $c$ which is also stable under the $T$-action determines a class $[Y]_T \in H^{2c}_T(X)$.

The set of $T$-fixed points on $X=G/P$ is identified with the coset space $W/W_P,$ where 
$W$ and $W_P$ are the Weyl groups of $G$ and $P$ respectively. 
Let us denote by $W^P$ the set of minimal length representatives for $W/W_P.$
We consider the case when $G$ is of one of the classical types B,C, and D, and $P$ is maximal. Then 
the set  $W^P$ will be later identified with certain set of strict partitions.  
For $\la\in W^P$, we denote $e_\la$ the corresponding $T$-fixed point on $X.$
Let $B$ be a Borel subgroup such that $T\subset B\subset P,$
and let $B^-$ be the opposite Borel subgroup.
A {\em Schubert variety} $\Omega_\lambda := \overline{B^- e_\lambda} \subset X$ is the closure of a $B^-$-orbit. 
The codimension of the variety $\Omega_\la$ is given by the  
length of $\la\in W^P$ denoted by $|\la|.$
The Schubert varieties are $T$-stable,  and
the {\em equivariant Schubert classes\/} $\sigma_\lambda:=[\Omega_\lambda]_T$ form an $S$-basis 
of $H^*_T(X)$, when $\lambda$ varies in $W^P$.

\subsection{Equivariant quantum cohomology} The equivariant quantum cohomology ring $\QH^*_T(X)$ is a graded $S[q]$-algebra, where $\deg t_i = 1$ and $q$ is a variable of degree $\deg q = c_1(T_X) \cap [C]$; here $c_1(T_X)$ is the first Chern class of the tangent bundle of $X$. (In this paper we use notation $\mathrm{deg}(\alpha)$ to indicate the complex degree of a homogeneous element $\alpha$ of $\QH^*_T(X).)$
The algebra $\QH^*_T(X)$ has an $S[q]$-basis given by Schubert classes $\sigma_\lambda$. The multiplication is given by the $3$-point, genus $0$ equivariant Gromov-Witten invariants $c_{\lambda, \mu}^{\nu,d}$, where $d \ge 0$ is a {\em degree} (a non-negative integer):
\[ \sigma_\lambda * \sigma_\mu = \sum_{\nu, d} c_{\lambda, \mu}^{\nu, d} q^d \sigma_\nu \/. \] This ring was defined by Givental and Kim \cite{GK,kim:oneqq}, and we refer to \cite{Mihalcea-EQSC, Mihalcea-EQChevalley} for more details about definitions, in the context of Grassmannians or homogeneous spaces. The coefficients $c_{\lambda,\mu}^{\nu,d}$ are homogeneous polynomials in $S$, and it was proved in \cite{mihalcea:positivity} that they can be written as positive sums in monomials of negative simple roots; see {\em loc.~cit.}~for precise details. This positivity generalizes the one in equivariant cohomology proved earlier by Graham \cite{graham:positivity}. The fact that equivariant quantum cohomology ring is a deformation of both the equivariant and quantum cohomology rings translates to the fact that if $d=0$ then $c_{\lambda,\mu}^{\nu,d}$ is the coefficient of $\sigma_\nu$ in the equivariant multiplication $\sigma_\lambda \cdot \sigma_\mu$, and if the degree of the polynomial $c_{\lambda,\mu}^{\nu,d}$ is zero (i.e. $c_{\lambda,\mu}^{\nu,d}$ is an integer) then the coefficient in question is the ordinary $3$-point, genus $0$ Gromov-Witten invariant which counts rational curves of degree $d$ intersecting general translates of varieties $\Omega_\lambda,\Omega_\mu$ and $\Omega_{w_0 \nu}$, where $w_0$ is the longest element in $W$.

We recall next a characterization theorem for the ring $\QH^*_T(X)$, proved in \cite{Mihalcea-EQChevalley}. Let $\sigma_{(1)} \in H^2_T(X)$ denote the {\em unique\/} Schubert class corresponding to the Schubert divisor (uniqueness follows because $P$ is maximal parabolic).

\begin{thm}\label{prop:chareqq}(\cite{Mihalcea-EQChevalley})
Let $(A,\star)$ be a graded commutative $S[q]$-algebra.
Assume that
\begin{enumerate}
\item $A$  has a $S[q]$-basis $\{s_\la\}_{\la\in W^P}$ such that 
$s_\la$ is homogeneous of degree $|\la|$;
\item The {\em equivariant quantum Chevalley rule\/} holds, i.e. 
$$s_{(1)} \star s_\lambda = \sum_{\mu, d\ge 0} c_{(1),\lambda}^{\mu,d} q^d s_\mu.$$
\end{enumerate}
Then $A$ is isomorphic to $QH^*_T(X)$
as a graded $S[q]$-algebra via the map defined by $s_\lambda\mapsto \sigma_\lambda.$
\end{thm}

The explicit formula for the equivariant quantum Chevalley rule was found in \cite{Mihalcea-EQChevalley}.
It states that there are no ``mixed" coefficients in $\sigma_{(1)} \star \sigma_\lambda$, i.e. all the coefficients are already founds in the equivariant, respectively quantum, specializations of the Chevalley rules. The equivariant coefficients appearing in this formula have been computed by Kostant and Kumar \cite{KK}; see also \cite{billey:polynomials}.~The case when $X$ is a maximal orthogonal Grassmannian - which is of main interest in this paper - was studied extensively in \cite{EYD,Ik}. For general homogeneous spaces $G/P$, the quantum Chevalley formula has been conjectured by Peterson \cite{peterson} and proved by Fulton and Woodward \cite{fulton.woodward} in its highest generality - see also \cite{buch.mihalcea:nbhds} for a different proof. Earlier results (including more general formulas for Grassmannians) were obtained by Bertram for type A Grassmannians \cite{bertram:quantum}, by Ciocan-Fontanine and Fomin, Gelfand and Postnikov for type A flag manifolds \cite{ciocan:partial,FGP} and by Kresch and Tamvakis for maximal isotropic Grassmannians \cite{KT-LG,KT-OG}; more recently a Pieri formula for submaximal isotropic Grassmannians was found in \cite{BKT:qiso}.

We recall next the definitions and the equivariant quantum Chevalley formulas for the maximal isotropic Grassmannians of types C and D, and how the relevant results in type B can be recovered from those of type D.

\subsection{Type C: the Lagrangian Grassmannian}\label{ss:C}The Lagrangian Grassmannian $\LG(n)$ is the manifold parametrizing dimension $n$ linear subspaces of $\C^{2n}$ which are isotropic with respect to a skew-symmetric, non-degenerate, bilinear form $\langle \cdot, \cdot \rangle_C $ on $\C^{2n}$. We fix an ordered basis $\e_n^*,\ldots,\e_1^*,\e_1,\ldots,\e_n$ of $\C^{2n}.$ The form is defined by 
\[\langle {\pmb e}_i, {\pmb e}_j\rangle_C  = \langle {\pmb e}_i^*, {\pmb e}_j^* \rangle_C = 0,\;
\langle {\pmb e}_i^*, {\pmb e}_j \rangle_C = \delta_{i,j} \/. \]
The symplectic group $\mathrm{Sp}_{2n}:=\mathrm{Sp}_{2n}(\C)$ acts transitively on $\LG(n)$.~In fact, $\LG(n)$ can be reinterpreted as the homogeneous space $\mathrm{Sp}_{2n}/P_0$ where $P_0 \subset \mathrm{Sp}_{2n}$ is a maximal parabolic subgroup which corresponds to the node $0$ of the Dynkin diagram of type $\mathrm{C}_{n}$:
\setlength{\unitlength}{0.4mm}
\begin{center}
  \begin{picture}(70,35)
  \thicklines
  \put(0,15){$\bullet$}
  \put(4,16.5){\line(1,0){12}}
  \put(4,18.5){\line(1,0){12}}
    \put(4.9,15){$>$}
  \multiput(15,15)(15,0){4}{
  \put(0,0){$\circ$}
  \put(4,2.4){\line(1,0){12}}}
  \put(75,15){$\circ$}
  \put(0,8){\tiny{$0$}}
  \put(15,8){\tiny{$1$}}
  \put(30,8){\tiny{$2$}}
  \put(72,8){\tiny{${n-1}$}}
  \put(50,8){{$\cdots$}}
  \end{picture}
\end{center}
~The Weyl group {$W=W_n$} of type $\mathrm{C}_n$ consists of {\em signed permutations\/} $w$ of the set $\{ \overline{n}, \ldots,\overline{1}, 1, \ldots, n\}$ which satisfy the property that $w(\overline{i}) = \overline{w(i)}$. Thus $w$ is determined by its values $w(1),\ldots,w(n)$. The minimal length representatives {$W_n^{P_0}$} correspond to (signed) {\em Grassmannian permutations\/}, which are defined by the property that $w(1) < w(2) <\cdots < w(n)$ in the ordering $\overline{n} < \cdots < \overline{1} < 1 < \cdots < n$. It follows that a signed Grassmannian permutation is completely determined by the subset of its barred values $w(1),w(2), \ldots,w(k)$, which in turn determines a {\em strict partition} $\lambda = (\lambda_1 > \cdots > \lambda_k)$ given by $\lambda_i = \overline{w(i)}$. Clearly $\lambda_1 \le n$ and $\lambda_k>0$. We denote this set by $\mathcal{SP}(n)$, and the set of {\em all} strict partitions  (i.e. with the requirement on $\lambda_1$ removed) by $\mathcal{SP}$.
  The identification ${W_n^{P_0}} \simeq \mathcal{SP}(n)$ is the same as that from \cite[\S 3]{DSP} or \cite[\S 4]{ikeda.naruse:K}, and we refer to any of these for more details, especially about the connection with the root theoretic description of $W$ and ${W_n^{P_0}}$.

\bigskip

\noindent{Example 1.} Let $n=6$. Then $\overline{3}\overline{6} 2 4 5 1$ is not a Grassmannian permutation, but $\overline{6}\overline{3} 1 2 4 5$ is. The latter determines the strict partition $\lambda= (6,3)$. 

\bigskip

The action of the maximal torus $T$ on $\C^{2n}$ determines a weight space decomposition $\C^{2n} = (\oplus_{i=1}^n \C {\pmb e}_i^*) \bigoplus (\oplus_{i=1}^n \C {\pmb e}_i)$ where $T$ acts by the character $t_i$ on $\C {\pmb e}_i^*$ (and $-t_i$ on $\C {\pmb e}_i$). 
Then $S=H^*_T(pt)$ equals the polynomial ring $\Z[t_1, \ldots, t_n]$, and in fact $t_i = c_1^T(\C {\pmb e}_i^*)$ (the equivariant first Chern class of a trivial line bundle). Let $F_i$ be the subspace spanned by the first 
$i$ vectors of the ordered basis. We have a complete flag \[ F_\bullet: F_1 \subset \cdots \subset F_n \subset F_{n+1} \subset \cdots \subset F_{2n} = \C^{2n}. \] Then $F_{i}$ is isotropic with respect to $\langle \cdot , \cdot \rangle_C$ for $1 \le i \le n$ and coisotropic for $n+1\le i\le 2n.$ In fact, if we denote
$V^\perp:=\{\u\in V\;|\;\langle \u,\v\rangle_C=0\;
\mbox{for all} \;\v\in V\},$ then $F_n^\perp=F_n$ and 
$F_{n+i}= F_{n-i}^\perp$ for $1\leq i\leq n$. The flag $F_\bullet$ is fixed by a Borel subgroup $B$, which in turn gives the opposite Borel subgroup $B^-$. Let $\lambda$ be a strict partition in $\mathcal{SP}(n)$ and $w_\lambda$ the corresponding Grassmannian element. Define $e_\lambda \in \LG(n)$ by $e_\lambda = \langle {\pmb e}_{w_\lambda(1)}^*, \ldots, {\pmb e}_{w_\lambda(n)}^* \rangle$, with the convention that ${\pmb e}_{\overline{i}}^* = {\pmb e}_i$. Then the Schubert variety $\Omega_\lambda = \overline{B^-  e_\lambda}$ can also be defined as  
\begin{equation} \Omega_\lambda = \{ V \in \LG(n)\;|\; \dim V \cap E^{\lambda_i} \ge i \quad(1\leq i\leq \ell(\la))\} \/.
\label{eq:Omega} 
\end{equation}
where $E^i=\langle \e_{i},\ldots,\e_n\rangle\:(1\leq i\leq n).$
With these conventions, the complex codimension of $\Omega_\lambda$ is $|\lambda| = \lambda_1 + \cdots + \lambda_{\ell(\lambda)}$.

We recall next the equivariant quantum Chevalley formula. In this case the quantum parameter $q$ has degree $n+1$.
For $\lambda \in \mathcal{SP}(n)$ such that $\la_1=n$ define 
$\la^{-}=(\la_2,\ldots,\la_{\ell(\la)})$. Then the following holds in $\QH^*_T(\LG(n))$:
\begin{equation}\label{E:eqqchevLG}\sigma_{(1)} *
\sigma_\la =\sum_{\mu\to \la,\;
\mu\in\mathcal{SP}(n) }2^{\ell(\lambda) - \ell(\mu) +1} \sigma_\mu
+c_{(1),\la}^\la
\sigma_\la+q\sigma_{\la^{-}},
\end{equation}
where 
$c_{(1),\la}^\la=2\sum_{i=1}^{\ell(\la)}t_{\la_i}$, 
and the last term is omitted if $\la^{-}$ does not exist.
This can be easily proved by using
the equivariant quantum Chevalley rule of \cite{Mihalcea-EQChevalley}
or \cite[Theorem 2.1]{LS}.

\subsection{Type D: the maximal orthogonal Grassmannian}\label{ss:D} We fix an ordered basis $\{{\pmb e}_{n+1}^*, \ldots, {\pmb e}_1^*, {\pmb e}_1, \ldots, {\pmb e}_{n+1} \}$ of $\C^{2n+2}$ and a symmetric, non-degenerate, bilinear form $\langle \cdot , \cdot \rangle_D$ which satisfies $\langle {\pmb e}_i, {\pmb e}_j \rangle_D = \langle {\pmb e}_i^*, {\pmb e}_j^* \rangle_D = 0$ and $\langle {\pmb e}_i^*, {\pmb e}_j \rangle_D = \delta_{i,j}$. 
Let $T$ be the maximal torus of the complex special orthogonal group $\SO_{2n+2}:=\SO_{2n+2} (\C)$
diagonalizing the basis. 
The maximal isotropic subspaces in $\C^{2n+2}$ have dimension $n+1$. Let $F_i$ denote the subspace spanned by the first 
$i$ vectors of the above basis; thus $F_{n+1}$ is maximal isotropic. The group $\SO_{2n+2}$ acts 
on the set of all maximal isotropic subspaces with two orbits, which correspond to the two connected components of this set. We denote by $\OG(n)$ the $\SO_{2n+2}$-orbit through $F_{n+1}.$
For a maximal isotropic subspace $V$, the condition that $V$ belongs to $\OG(n)$
is equivalent to the following {\em even parity condition\/}:
\begin{equation}
\mathrm{codim}_{F_{n+1}} (V \cap F_{n+1})\equiv 0\;\mod 2.\label{eq:parity}
\end{equation}
The other $\SO_{2n+2}$-orbit corresponds to the odd parity condition.

%For an isotropic subspace of $\C^{2n+2}$, the condition that it is maximal
%among the isotropic subspaces is equivalent to be of dimension $(n+1).$
%Let $F_i$ denote the subspace spanned by the first 
%$i$ vectors of the above basis. Then $F_{n+1}$ is clearly a maximal isotropic subspace. 
%Obviously, the group $\SO_{2n+2}$ acts 
%on the set of all maximal isotropic subspaces. 
%Let us denote by $\OG(n)$ the $\SO_{2n+2}$-orbit through $F_{n+1}.$
%For a maximal isotropic subspace $V$, the condition that $V$ belongs to $\OG(n)$
%is equivalent to 
%the following {\em even parity condition\/}:
%\begin{equation}
%\mathrm{codim}_{F_{n+1}} (V \cap F_{n+1})\equiv 0\;\mod 2.\label{eq:parity}
%\end{equation}
%There is another $\SO_{2n+2}$-orbit corresponding to the odd parity condition.
%Note that 
%the two orbits are also a two connected components of the set of all maximal isotropic subspaces of $\C^{2n+2}.$

Each of these two components is a homogeneous space. In fact, 
$\OG(n)$ can be identified with $\SO_{2n+2}/P_{\hat{1}}$, while the other component with  
 $\SO_{2n+2}/P_1$,
 where $P_1$ and $P_{\hat{1}}$ are the maximal parabolic subgroups determined by the ``forked" nodes $1, \hat{1}$ of the Dynkin diagram of type $\mathrm{D}_{n+1}$
 \footnote{We use following convention for the simple roots, which 
is stable with respect to the natural embedding of a Dynkin diagram to ones of higher ranks:  
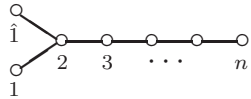
 \[\alpha_{\hat{1}}=t_2 + t_1,\quad \alpha_1= t_2 - t_1,\quad \alpha_2=t_3-t_2,  \ldots,\quad \alpha_n= t_{n+1} - t_{n}. \] }:
\setlength{\unitlength}{0.4mm}
\begin{center}
  \begin{picture}(70,35)
  \thicklines
  \put(0,0){
  \put(0,25){$\circ$}
  \put(0,5){$\circ$}
  \put(4,26){\line(3,-2){12}}
  \put(4,8.5){\line(3,2){12}}
  \multiput(15,15)(15,0){4}{
  \put(0,0){$\circ$}
  \put(4,2.4){\line(1,0){12}}}
  \put(75,15){$\circ$}
  \put(0,17){\tiny{$\hat{1}$}}
  \put(0,-1){\tiny{$1$}}
  \put(16,8){\tiny{$2$}}
  \put(31,8){\tiny{$3$}}
  \put(75,8){\tiny{${n}$}}
  \put(46,8){$\cdots$}}
  \end{picture}
\end{center}

{
\noindent{Remark 3.}  }
Let $\iota$ denotes the automorphism of $G$ corresponding to the automorphism of the Dynkin diagram switching nodes $1$ and $\hat{1}$.
The automorphism $\iota$ induces an isomorphism of algebraic varieties $\varphi: G/P_1 \to G/P_{\hat{1}}$, which
is {not} equivariant with respect to the natural $T$-actions, 
but it satisfies $\varphi(t \cdot x) = \iota(t) \cdot \varphi(x)\;(t\in T)$. 

\bigskip

 %(regarded as characters of $T$ written additively).
%\begin{footnote}{Note that this choice of simple roots  $\mathrm{D}_n \subset \mathrm{D}_{n+1}$.}\end{footnote} 
%Then $P_{\hat{1}}$ is the maximal parabolic group associated with the simple root $\alpha_{\hat{1}}$. Using the conventions from \cite[\S 10]{DSP}, the variety $G/P_{\hat{1}}$ consists of the maximal isotropic subspaces $V \subset \C^{2n +2}$ such that the codimension of $V \cap F_{n+1}$ in $F_{n+1}$ . We refer to this as the {\em maximal orthogonal Grassmannian}, and we denote it by $\OG(n)$. 
 
The Weyl group $W'$ of type $\mathrm{D}_{n+1}$ is the subgroup of the type $\mathrm{C}_{n+1}$ Weyl group from the previous section which consists of signed permutations of $\{ \overline{n+1}, \ldots, \overline{1}, 1, \ldots, n+1 \}$ {\em with even number of sign changes}. This identifies the simple reflection $s_{\hat{1}}:=s_{t_2 + t_1}$ with the permutation $\overline{2}\overline{1} 3\cdots (n+1)$ and the simple reflection $s_i:=s_{t_{i+1} - t_i}$ with the element $s_i:= (i, i+1)(\overline{i},\overline{i+1})$. The set of minimal length representatives for $W'/W'_{P_{\hat{1}}}$ coincides with the set of those Grassmannian permutations from type $\mathrm{C}_{n+1}$ which are included in $W'$. For such a representative $w$ for $W'/W'_{P_{\hat{1}}}$ define the strict partition 
$\lambda=(\lambda_{1}>\cdots>\lambda_{r}\geq 0)$ by setting $ \lambda_{i}=\overline{w(i)}-1$ where $i$ varies over the set $1, \ldots, r$ such that $w(i)< 1$.
Note that $r$ must be even and $\lambda_{r}$ can be zero
this time. This correspondence gives a
bijection with the set
$\mathcal{SP}(n)$ and we denote by $w'_{\lambda}$
the representative corresponding
to $\lambda\in \mathcal{SP}(n)$.
As before, $\ell(w'_{\lambda})=|\lambda|$ where 
$\ell(w)$ denotes the length of $w$ as an element of the Weyl group of type $D_{n+1}$. We refer again to \cite{DSP} or \cite{ikeda.naruse:K} for more details.

\bigskip

\noindent{Example 2.} Let $\lambda=(4,2,1).$
Then the corresponding Grassmannian elements are given by
$w_{\lambda}=\bar{4}\bar{2}\bar{1}3
=s_{0}s_{1}s_{0}s_{3}s_{2}s_{1}s_{0}$ in type C
and $w'_{\lambda}=\bar{5}\bar{3}\bar{2}\bar{1}4
=s_{\hat{1}}s_{2}s_{1}s_{4}s_{3}s_{2}s_{\hat{1}}$ in type D.

\bigskip

As in type C, we can choose a maximal torus $T \subset \SO_{2n+2}$ which acts on $\C^{2n+2}$ with ${\pmb e}_i, {\pmb e}_i^*$ as eigenvectors, and such that the weight of $T$ on $\C {\pmb e}_i^*$ and $\C {\pmb e}_i$ is $t_i$, respectively $-t_i$. Then $S=H^*_T(pt) = \Z[t_1, \ldots, t_{n+1}]$ and geometrically $t_i = c_1^T(\C {\pmb e}_i^*)$. If $\lambda$ is a strict partition in $\mathcal{SP}(n)$ let $e_\lambda = \langle {\pmb e}_{w'_\lambda(1)}^*, \ldots, {\pmb e}_{w'_\lambda(n+1)}^* \rangle$, using again the convention that ${\pmb e}_{\overline{i}}^* = {\pmb e}_i$. Let $B^{-}$ be the Borel subgroup opposite to the stabilizer of $F_\bullet.$
Then the Schubert variety is defined by $\Omega_\lambda = \overline{B^-  e_\lambda}$, and it also equals 
\[\Omega_\lambda = \{ V \in \OG(n)\;|\; \dim V \cap E^{\lambda_i} \ge i \quad(1\leq i\leq \ell(\la))\},\]
where $E^i=\langle \e_{i+1},\ldots,\e_{n+1}\rangle \; (1\leq i\leq n).$
The Schubert variety determines an equivariant Schubert class $\sigma_\lambda \in H^{2 |\lambda|}_T(\OG(n))$. In this case, the quantum parameter has degree $\deg q = 2n$. If $\la$ is a strict partition such that 
$\la_1=n$ and $\la_2=n-1$ then define $\lambda^- =( \lambda_3, ..., \lambda_r)$. Otherwise we will say that $\lambda^-$ does not exist. With these notations, the equivariant quantum Chevalley formula is given by \begin{equation}\label{E:eqqchevOG} 
\sigma_{(1)} * \sigma_\la=\sum_{\mu\to \la,\;
\mu\in\mathcal{SP}(n) }\sigma_\mu
+c_{(1),\la}^\la
\sigma_\la+q\sigma_{\la^{-}},
\end{equation}
where $c_{(1),\la}^\la=\sum_{i=1}^{r}t_{\la_{i}+1}$,
and the last term is omitted if $\la^{-}$ does not exist.
This can be easily proved by using
the equivariant quantum Chevalley rule of \cite{Mihalcea-EQChevalley}.

\subsection{The relation between maximal orthogonal Grassmannians in types B and D}\label{ssec:BD}
Let $\langle \cdot , \cdot \rangle_B$ be a non-degenerate, symmetric bilinear form on $\C^{2n+1}$. The maximal dimension of a subspace $V \subset \C^{2n+1}$ which is isotropic with respect to $\langle \cdot , \cdot \rangle_B$ equals $n$. The variety $\OG(n, \C^{2n+1})$ parametrizing maximal isotropic subspaces in $\C^{2n+1}$ is called again the maximal orthogonal Grassmannian. It is a homogeneous variety $\SO_{2n+1}/P$, where $P$ is the maximal parabolic group determined by the end root corresponding to the double edge of the Dynkin diagram of type B. Using that the Weyl groups of type B and C coincide, one shows that the Schubert varieties in $\OG(n, \C^{2n+1})$ are indexed again by strict partitions in $\mathcal{SP}(n)$. 

We turn to the relation between the equivariant quantum cohomology rings of $\OG(n, \C^{2n+1})$ and $\OG(n)$. Fix a basis $\{ {\pmb e}_{n+1}^*, \ldots, {\pmb e}_2^*, {\pmb e}_0, {\pmb e}_2, \ldots, {\pmb e}_{n+1} \}$ for $\C^{2n+1}$ such that 
\[ \langle {\pmb e}_i, {\pmb e}_j \rangle_B = \langle {\pmb e}_i^*, {\pmb e}_j^* \rangle_B = \langle {\pmb e}_0
, {\pmb e}_i^* \rangle_B =\langle {\pmb e}_0, {\pmb e}_j \rangle_B= 0 \; (2\leq i\leq n),\]
\[ \langle {\pmb e}_i^*, {\pmb e}_j \rangle_B = \delta_{i,j}\; (2\leq i,j\leq n)\quad \langle {\pmb e}_0, {\pmb e}_0 \rangle_B =1 \/. \]
Consider the injective linear map  $\C^{2n+1} \to \C^{2n+2}$ which identifies the vectors ${\pmb e}_i^*,\;{\pmb e}_i\;(2\leq i\leq n)$ from $\C^{2n+1}$ with those from $\C^{2n+2}$ from the previous section, and sends ${\pmb e}_0 \in \C^{2n+1}$ to $\frac{1}{\sqrt{2}} ({\pmb e}_1^* + {\pmb e}_1) \in \C^{2n+2}$. Under this injection, $\C^{2n+1}$ can be identified with the space $(\C ({\pmb e}_1^* - {\pmb e}_1))^\perp \in \C^{2n+2}$ orthogonal to $\C ({\pmb e}_1^* - {\pmb e}_1)$ with respect to the symmetric form $\langle \cdot , \cdot \rangle_D$ on $\C^{2n+2}$ defined in the section \ref{ss:D} above. This identifies $\SO_{2n+1}$ with a subgroup of $\SO_{2n+2}$, and the symmetric form $\langle \cdot , \cdot \rangle_D$ from type D restricts to $\langle \cdot , \cdot \rangle_B$ on $\C^{2n+1}$. 
Note that $\C\e_1$ and $\C\e_1^*$ are 
 exactly the two maximal isotropic subspaces in $\langle{\pmb e}_1^*, {\pmb e}_1\rangle\cong \C^2$. 
It follows that if $V\in  \OG(n,\C^{2n+1})$, then  
$V+ \C\e_1$ and $V+ \C\e_1^*$ 
are the only maximal isotropic subspaces in $\C^{2n+2}$
containing $V.$ Only one of them
satisfies 
the parity condition (\ref{eq:parity}) and gives an element of $\OG(n)$.
This correspondence gives an isomorphism $\eta:\OG(n,\C^{2n+1}) \to \OG(n)$.
This isomorphism is equivariant with respect to the inclusion of maximal tori of $\SO_{2n+1}$ and $\SO_{2n+2}$ given by $\diag(\xi_{n+1}^{-1}, \ldots, \xi_2^{-1}, 1, \xi_2, \ldots, \xi_{n+1}) 
\mapsto \diag(\xi_{n+1}^{-1}, \ldots, \xi_2^{-1}, 1,1, \xi_2, \ldots, \xi_{n+1}) $. 
Let $\Omega_\la'$ be the Schubert variety for $\OG(n,\C^{2n+1})$ associated with $\la\in \mathcal{SP}(n).$ 
Then we have 
$\eta^{-1}(\Omega_\la)=\Omega_\la'$.
Combining all of the above implies that there exists a $\Z[t_1, \ldots, t_{n+1}][q]$-algebra isomorphism \[ \eta^*: \QH^*_T(\OG(n))/\langle t_1 \rangle  \simeq \QH^*_{{T'}}(\OG(n, \C^{2n+1})) \] so that $\eta^*[\Omega_\lambda]_T = [\Omega_{\lambda}']_{T'}$, where $T':=T\cap \SO_{2n+1}$.
Note that $H_{T'}^*(pt)=\Z[t_2,\ldots,t_{n+1}]$.
This explicit isomorphism allows us to ignore the type B orthogonal Grassmannian from now on, and instead focus on the one of type D.

\setcounter{equation}{1}
\section{Presentations and Giambelli formulas for the equivariant quantum ring of the orthogonal Grassmannian}\label{sec:OG}

The goal of this section is to prove the main results of this paper in the case of the orthogonal Grassmannian $\OG(n)$. Our strategy is to work with both the rings given by the ``abstract" variables and by the factorial $P$-Schur functions at the same time, and show that each ring has a natural $S[q]$-basis for which the equivariant quantum Chevalley formula (\ref{E:eqqchevOG}) holds.

\subsection{Equivariant Quantum Chevalley rule in $\mathcal{A}_n$}
In this section, we set $t_i=0$ for $i>n+1.$
Define the ring
\begin{equation*}
\mathcal{A}_n:= S[q][P_1(x|t), \ldots, P_{2n}(x|t)]/I_n^{(q)},\;\;
I_n^{(q)}:=\langle P_{n+1}(x|t), \ldots, P_{2n-1}(x|t), P_{2n}(x|t) + (-1)^n q \rangle.
\end{equation*}
We will show that the images of $P$-Schur functions %$P_\lambda(x|t) \in \mathcal{A}_n$ 
satisfy the equivariant quantum Chevalley rule %(\ref{E:eqqchevOG}) 
for $\OG(n)$.

\begin{thm}\label{prop:EQCh}
Let $\la \in \mathcal{SP}(n)$ be a strict partition with at most $n$ parts. Then
$$
P_{(1)}(x|t){P_\la(x|t)}
 \equiv 
 \sum_{\mu\to \la,\;\mu\in \mathcal{SP}(n)
 }{P_\mu(x|t)}
+
c_{(1),\la}^\la(t) {P_\la(x|t)} +q\,{P_{\la^{-}}(x|t)}\quad \mod I_n^{(q)},
 $$
 where 
 the last term omitted unless $\la_1=n$ and $\la_2=n-1.$ 
\end{thm}

\noindent{Proof.} By the Chevalley rule (\ref{eq:Chev}), we have 
$$
P_{(1)}(x|t){P_\la(x|t)}=
\sum_{\mu\to \la,\;\mu \in \mathcal{SP}(n)}{P_\mu(x|t)}+c_{(1),\la}^\la(t) {P_\la(x|t)}
+{P_{\overline{\la}}(x|t)},
 $$
where if $\la_1=n$ then $\overline{\la}$ is given by adding one box
 to the first row of $\la$
 otherwise the last term is omitted. %if $\la_1<n.$
 Let $\nu\in \mathcal{SP}(n-2)$ and denote by $(n+1, k, \nu)$ the strict partition with parts $n+1, k$ followed by the parts of $\nu$ (thus $k > \nu_1$). It is enough to show 
$$P_{(n+1,k,\nu)}(x|t)\equiv \begin{cases}
q P_{\nu}(x|t)\quad
\mod I_n^{(q)}&  \mbox{if}\;k=n-1\\
0 \quad
\mod I_n^{(q)}& \mbox{if}\;k<n-1.
\end{cases}$$
If $k<n-1$ 
then the equation (\ref{eq:rec}) and induction on $\ell(\lambda)$ implies that
$P_{(n+1,k,\nu)}(x|t)$ belongs to $I_n^{(q)}$
because it can be 
expressed as an $S$-linear combination of 
$P_r(x|t)P_s(x|t)$ with $n+1\leq r<2n.$ Let now $\nu=\emptyset$ and $k=n-1$. Then by Proposition \ref{prop:quadratic}
$$
P_{(n+1,n-1)}(x|t)=P_{n+1}(x|t)P_{n-1}(x|t)
+\sum_{(r,s) \in \mathcal{I}_{n+1,n-1}} g_{n+1,n-1}^{r,s}(t)P_{r}(x|t)P_s(x|t)\/. $$
From the definition of $\mathcal{I}_{n+1,n-1}$ it follows that the right hand side is equivalent to 
$(-1)^{n-1}P_{2n}(x|t)
\equiv q$ modulo $I_n^{(q)}$. For general $\nu$ and $k=n-1$, invoking again (\ref{eq:rec}), 
we can express $P_{(n+1,n-1,\nu)}(x|t)$ 
as $$
P_{(n+1,n-1,\nu)}(x|t)
=P_{(n+1,n-1)}(x|t)P_{\nu}(x|t)+\mbox{other terms}.
$$
All the other terms have factor $P_{(n+1,k)}(x|t)$
with $k<n-1$, therefore are congruent to zero
by the same reason as in the first case.
Since we have $P_{(n+1,n-1)}(x|t)\equiv q\;\mod I_n^{(q)}$
we have completed the proof.
\qed

\subsection{Presentation of $\mathcal{A}_n$ as a quotient ring}\label{ssec:twomodels} 
Next we will obtain a presentation for 
$\mathcal{A}_n$ as a quotient ring of a polynomial ring. 
Let $X_1,\ldots,X_{2n}$ be indeterminates.
Consider the $S$-algebras homomorphism
\[S[X_1,\ldots,X_{2n}]\rightarrow
 S[P_1(x|t),\ldots,P_{2n}(x|t)]
\] sending $X_i$ to $P_i(x|t)$, for $1 \le i \le 2n$.
The quadratic relation (\ref{eq:quadratic}) implies that 
the kernel of the map contains the quadratic polynomials
$R_i$  in $S[X_1, \ldots, X_{2n}]$ defined by \begin{equation}
R_i= X_i^2+\sum_{(r,s)}g_{i,i}^{r,s}(t)
X_rX_s,\label{eq:R_i}
\end{equation}
where the sum is over the pairs $(r,s) \in \mathcal{I}_{(i,i)}$ defined in equation (\ref{eq:rs}), and coefficients $g_{i,i}^{r,s} \in S$ are defined in (\ref{E:g}) above. 

\begin{lem}\label{lemma:univiso} 
The morphism  
\[S[X_1,\ldots,X_{2n}] /\langle R_1,\ldots,R_n\rangle\rightarrow
 S[P_1(x|t),\ldots,P_{2n}(x|t)]\]
  sending $X_i$ to $P_i(x|t)$ for $1 \le i \le 2n$ is an isomorphism of $S$-algebras. \end{lem}

\noindent{Proof.} It suffices to show injectivity. Using the quadratic relation $R_i$ and induction on $i$, one can write $X_{2i}\;(1\leq k\leq n)$ as a polynomial in $X_1,X_3,\ldots,{X_{2k-1}}$
with coefficients in $S$. 
Then $S[X_1,\ldots,X_{2n}]/\langle R_1,\ldots,R_n\rangle$ is
generated as an $S$-module by the set of the residue classes of the monomials 
\begin{equation}
{X_1}^{m_1}{X_3}^{m_3}\cdots 
{X}_{2n-1}^{m_{2n-1}}\quad( m_1,m_3,\ldots,m_{2n-1}\geq 0)\/.\label{eq:monomials}
\end{equation}
Since 
$P_1(x),P_3(x),\ldots,P_{2n-1}(x)$
are algebraically independent over $\Z$ \cite[pag. 252]{Mac} we deduce that deformed functions 
$ P_1(x|t),P_3(x|t),\ldots,P_{2n-1}(x|t)$
are algebraically independent over $S$.
Then the image of the monomials from (\ref{eq:monomials}) is linearly independent
over $S$. Hence the map is injective.
\qed

\bigskip

Let us denote the quotient ring 
\begin{equation}
{\mathcal{B}}_n:= S[q][X_1,X_2,\cdots,X_{2n}]/\langle 
X_{n+1},\ldots,X_{2n-1},X_{2n}+(-1)^nq,
R_1,\ldots,R_n\rangle.
\end{equation}
\begin{prop}\label{prop:A=B}
We have an isomorphism
of graded $S[q]$-algebras \[
\mathcal{B}_n \to \mathcal{A}_n
\]
sending $X_i$ to $P_i(x|t)\;(1\leq i\leq 2n).$
\end{prop}
\noindent{Proof.} Clear from the definitions of $\mathcal{A}_n,\mathcal{B}_n$ and 
Lemma
\ref{lemma:univiso}.
\qed

\subsection{Freeness of $\mathcal{B}_n$ over $S[q]$}\label{ssec:free} Next we will prove that $\mathcal{B}_n$ is
a free module over $S[q]$. The techniques are very similar to those from \cite[\S 8.2]{Ik}, where the (non-quantum) equivariant cohomology is considered. 

Let $R_i^{(q)} \in S[X_1,\ldots,X_n]$ be 
obtained from $R_i$  by substituting 
 $X_{n+1}=\ldots =X_{2n-1}=0$
and $X_{2n} = (-1)^{n-1}q$. For example, if $n=2$, then $S=\Z[t_1,t_2,t_3]$ and 
\begin{eqnarray*}
R_1&=& X_1^2-X_2+(t_1+t_2)X_1;\\
R_2& =& X_2^2 - 2X_3X_1 + X_4 + 2t_3X_2X_1- (t_1 + t_2 + t_3)X_3 + (t_1t_2 + t_1t_3 + t_2t_3 + t_3^2)X_2.
\end{eqnarray*}
Thus we have 
\begin{equation*}
R_1^{(q)}=R_1,\quad
R_2^{(q)} =X_2^2
 - q+2t_3 X_2X_1+(t_1t_2+t_1t_3+t_2t_3+t_3^2)X_2.
\end{equation*}
With this notation we have the following 
identification \begin{equation}\label{E:eq} \mathcal{B}_n = S[q][X_1, \ldots, X_n]/\langle R_1^{(q)}, \ldots, R_n^{(q)} \rangle \/. \end{equation}

\noindent For $\la=(\la_1,\ldots,\la_{\ell(\la)})\in \mathcal{SP}(n) $, set
 $X^\la=X_{\la_1}\cdots X_{\la_{\ell(\la)}}$.

\begin{lem}\label{lem:free} The ring 
$\mathcal{B}_n$ is 
a free $S[q]$-module  
with a basis 
\begin{equation}\{X^\la\;|\;\la\in\mathcal{SP}(n)\}.
\label{eq:dist-monomials}
\end{equation}
\end{lem}
\noindent{Proof.} Let $X^{\bf m} = X_1^{m_1}X_2^{m_2}\cdots
X_n^{m_n}$ be an arbitrary monomial
in $X_1,\ldots,X_n$. Using the ``quantized" quadratic relations $R_i^{(q)}$, we will write this as an $S[q]$-linear combination of the square-free monomials from (\ref{eq:dist-monomials}). The argument is almost identical to that from \cite[Lemma 8.4]{Ik}.
Consider the graded reverse lexicographic 
order with $X_1<\cdots<X_n$. The key property of this ordering is that \[ R_i^{(q)} = X_i^2 +  \mathrm{lower~ order~ terms} \/. \] If $X^{\bf m}$ has 
an exponent $m_i\geq 2$ we replace $X_i^2$ by an $S[q]$-
linear combination of monomials which are all smaller than
$X_i^2$. Repeating this yields the claimed square-free linear combination. Thus the set (\ref{eq:dist-monomials}) spans
$\mathcal{B}_n$ over $S[q]$.

We prove next that the set  (\ref{eq:dist-monomials}) is linearly 
independent over $S[q]$ by using 
the theory of Gr\"obner bases. 
Let $K$ be the field of fraction of $S$. By Buchberger criterion (see e.g. Exercise 15.20 and more generally \cite[Ch. 15]{eisenbud:comalg}) $R_1^{(q)}, \ldots, R_n^{(q)}$ is a Gr\"obner basis for the ideal $\langle R_1^{(q)}, \ldots, R_n^{(q)} \rangle \subset K[q, X_1, \ldots, X_n]$
for any monomial ordering such that the initial term of $R_i^{(q)}$ is $X_i^2$. One such ordering is the graded reverse lexicographic order with $q<X_1<\cdots<X_n$ (where we temporarily declare that $\deg q =1$). Notice that  
$\{q^dX^\la\;|\;d\geq 0,\;
\la\in \mathcal{SP}(n)\}$
consists exactly of the  
monomials in $q,X_1,\ldots,X_n$ 
which are not divisible by $X_1^2,X_2^2,\ldots,X_n^2$.
It follows that this set is linearly independent
over $K$ in the quotient ring
$K[q,X_1,\ldots,X_n]/\langle R_1^{(q)}, \ldots, R_n^{(q)} \rangle$.
From this fact we see that (\ref{eq:dist-monomials}) 
is linearly independent over $S[q]$ as well.
\qed

\subsection{The Pfaffian polynomials $X_\la$}\label{ssec:X_la}
We give another $S[q]$-basis of 
$\mathcal{B}_n$, which 
will be finally identified with the Schubert basis. 
In analogy with the Pfaffian formula  (\ref{eq:Pfaff-P}) for $P_\lambda(x|t)$,  define $X_\la\in S[q][X_1,\ldots,X_n]$ by  
\begin{equation}\label{E:Xla}
X_\la=\mathrm{Pf} (X_{\la_i,\la_j})_{1\leq i<j\leq r}; \quad X_{k,l}=X_kX_l+\sum_{(r,s) \in \mathcal{I}_{k,l} }g_{k,l}^{r,s}(t)
X_rX_s 
\end{equation}
where $\mathcal{I}_{r,s}$ is defined in (\ref{eq:rs}) above, and we make the substitutions \[ X_{n+1} = \ldots = X_{2n-1} =0; \quad X_{2n} = (-1)^{n-1}q \/. \] 

\medskip

\noindent{Example 3.} If $n=2$,  
we have 
$P_{2,1}(x|t)
=P_2(x|t)P_1(x|t)
-P_3(x|t)-(t_1+t_3)P_2(x|t),
$
thus $X_{2,1}=X_2X_1-(t_1+t_3)X_2$. 

\bigskip

The definition of $X_\lambda$ and the Pfaffian formula (\ref{eq:Pfaff-P}) imply that the isomorphism $\mathcal{B}_n \to \mathcal{A}_n$ from Proposition \ref{prop:A=B} sends $X_\lambda$ to $P_\lambda(x|t)\;\mod I_n^{(q)}$, for any $\lambda \in \mathcal{SP}(n)$. 

\begin{prop}\label{T:Ybasis} The polynomials $X_\lambda$, when $\lambda$ varies in $\mathcal{SP}(n)$ form a $S[q]$-basis for 
$\mathcal{B}_n= S[q][X_1,\ldots,X_n]/\langle R_1^{(q)}, \ldots, R_n^{(q)} \rangle$. \end{prop}

\noindent{Proof.} Consider the expansion $
X_\la=
\sum_{\mu\in \mathcal{SP}(n)}c_{\la\mu} X^\mu$,
where $c_{\lambda \mu}$ is a homogeneous polynomial in $S$ of degree $|\la| - |\mu|$. A standard argument using the Pfaffian recurrence from (\ref{eq:rec}) and induction on the number of parts of $\lambda$ shows that $c_{\la\la}=1$ and $c_{\la\mu}=0$ unless
$\mu\leq \la$; see \cite[Lemma 8.5]{Ik} for more details. Then the matrix 
$(c_{\la\mu})$ is invertible, and the result follows from Lemma \ref{lem:free}.
\qed

\subsection{Proof of the main Theorem in type D}
\begin{thm}\label{thm:OG} (i) (the factorial $P$-Schur version) 
There is an isomorphism of 
graded $S[q]$-algebras \[
S[q][P_1(x|t),\ldots,P_{2n}(x|t)]/I_n^{(q)}
\longrightarrow \QH_T^*(\OG(n)),
\]
where 
$I_n^{(q)}$ is the ideal generated by 
$P_{n+1}(x|t),\ldots,P_{2n-1}(x|t), P_{2n}(x|t)+(-1)^n q$.
Moreover, the image of $P_\la(x|t)\;(\la \in \mathcal{SP}(n))$
is the Schubert class $\sigma_\la.$ 

(ii) (the indeterminate version) There is an isomorphism of 
graded $S[q]$-algebras
$$
S[q][X_1,X_2,\ldots,X_n]/\langle 
R_1^{(q)},R_2^{(q)},\ldots,R_n^{(q)}\rangle
\longrightarrow \QH_T^*(\mathrm{OG}(n))
$$
which sends $X_\la\;(\la \in \mathcal{SP}(n))$
to the Schubert class $\sigma_\la$. 
\end{thm}
\noindent{Proof.} By Proposition \ref{prop:A=B} there is an isomorphism of graded $S[q]$-algebras $\mathcal{B}_n \simeq \mathcal{A}_n$, which sends $X_\lambda,\;\la\in\mathcal{SP}(n)$ to $P_\lambda(x|t) \;\mod I_n^{(q)}$. The theorem follows by combining the facts that the equivariant quantum Chevalley formula holds in $\mathcal{A}_n$ by Theorem \ref{prop:EQCh}, $\mathcal{B}_n$ is free over $S[q]$ with the ``correct" Schubert basis (Proposition \ref{T:Ybasis}), and by applying the characterization theorem \ref{prop:chareqq}. 
\qed

\setcounter{equation}{1}
\section{Equivariant quantum cohomology of the Lagrangian Grassmannian}\label{sec:LG}
 In this section we study the presentation and the Giambelli formulas for the equivariant quantum ring of $\LG(n)$. The arguments are similar to those from type D,
although some technical difficulties arise. 
\subsection{The Chevalley rule} 
In this section, we put $t_i=0$ for $i>n.$
Recall that $S=\Z[t_1, \ldots, t_n]$.
Let $\mathcal{A}_n$ be the ring defined by the left hand side of Theorem \ref{thm:main} (b). 
\begin{thm}[Equivariant Quantum Pieri-Chevalley rule] \label{T:EQCh-C}
Let $\la$ be a strict partition in $\mathcal{SP}(n)$. Then
\[ Q_{(1)}(x|t){Q_\la(x|t)}
 \equiv 
 \sum_{\mu\to \la,\;
 \mu\in\mathcal{SP}(n)}2^{\ell(\lambda) - \ell(\mu) +1} {Q_\mu(x|t)}
+  2(\sum_{i=1}^r t_{\la_i+1}){Q_\la(x|t)} +q{Q_{\la^{-}}(x|t)} \quad \mod J_n^{(q)},\]
where the last term omitted unless $\la_1=n$. 
\end{thm}

\noindent{Proof.} The idea of the proof is similar to that of Theorem \ref{prop:EQCh}, but the actual details are slightly different. First, given the Chevalley rule from Proposition (\ref{eq:Chev}) for factorial $Q$-Schur functions, it suffices to show that for any partition $\nu \in \mathcal{SP}(n-1)$
\begin{equation}\label{E:claim}
2Q_{n+1,\nu}(x|t)
\equiv qQ_{\nu}(x|t) ~\mod J_n^{(q)} \/. 
\end{equation}
 Using the equation (\ref{E:Qrec}) and the fact that in this case $t_{n+1} = 0$ we obtain that $Q_{n+1,\nu_i}(x|t)\equiv 
Q_{n+1}(x|t)Q_{\nu_i}(x|t) ~\mod J_n^{(q)}$.

Assume first $\ell:=\ell(\nu)$ is odd. Then by (\ref{eq:rec}) 
\begin{equation}\begin{split}
Q_{n+1,\nu}(x|t)
&=
\sum_{i=1}^\ell(-1)^{i-1}Q_{n+1,\nu_i}(x|t)Q_{\nu_1,\ldots,\widehat{\nu_i},\ldots,\nu_\ell}(x|t)  \\ 
\equiv Q_{n+1}(x|t)& \sum_{i=1}^\ell (-1)^{i-1} Q_{\nu_i} (x|t) Q_{\nu_1,\ldots,\widehat{\nu_i},\ldots,\nu_\ell}(x|t) = Q_{n+1}(x|t) Q_\nu(x|t)
\end{split}\end{equation} where the last equality follows from (\ref{eq:rec-odd}). This implies (\ref{E:claim}).

Let now $\ell:=\ell(\nu)$ be even. Then the length of $(n+1, \nu)$ is odd and we apply (\ref{eq:rec-odd}): \begin{eqnarray*}
Q_{n+1,\nu}(x|t)
&=&
Q_{n+1}(x|t)Q_{\nu}(x|t)
+\sum_{i=1}^\ell(-1)^i
Q_{\nu_i}(x|t)Q_{n+1,\nu_1,\ldots,\widehat{\nu_i},\ldots,\nu_\ell}(x|t)\\
&\equiv &
Q_{n+1}(x|t)Q_{\nu}(x|t)
+Q_{n+1}(x|t)\sum_{i=1}^\ell(-1)^i
Q_{\nu_i}(x|t)Q_{\nu_1,\ldots,\widehat{\nu_i},\ldots,\nu_\ell}(x|t)\end{eqnarray*}
Now note that $\sum_{i=1}^\ell(-1)^i
Q_{\nu_i}(x|t)Q_{\nu_1,\ldots,\widehat{\nu_i},\ldots,\nu_\ell}(x|t)=0$. Indeed, by applying (\ref{eq:rec-odd}) 
to $Q_{\nu_1,\ldots,\widehat{\nu_i},\ldots,\nu_r}(x|t)$ we see that the coefficient of $ Q_{\nu_1,\ldots,\widehat{\nu_i},\ldots ,\widehat{\nu_j},\ldots,\nu_r}(x|t)$ in this sum
appears twice with opposite signs. This gives (\ref{E:claim}) in this case and finishes the proof. \qed

\subsection{Freeness and the Pfaffian formula} In this section we denote by $R_i \in S[X_1, \ldots, X_{2n}]$ the quadratic relations obtained from the analogous relations (\ref{E:Qqvanishing}) for the functions $Q_i(x|t)$: \[ R_i = X_i^2 +  \sum_{(r,s) \in \mathcal{I}_{i,i}} f_{i,i}^{r,s} X_r X_s; \quad 1 \le i \le n \/.\] The ``quantized" quadratic relations $R_i^{(q)}$ are obtained from $R_i$ using the substitutions \begin{equation}\label{E:LGsub}2X_{n+1} =q, X_{n+2} = \cdots = X_{2n} = 0 \/. \end{equation} Note that each monomial in $R_i$ containing $X_{n+1}$ has coefficient $2$, therefore $R_i^{(q)} \in S[q][X_1, \ldots, X_n]$.

\bigskip

\noindent{Example 4.} When $n=2$, we have \begin{eqnarray*}
R_1&=& R_1^{(q)} = X_1^2-2X_2-2t_1X_1,\\
R_2&=&X_2^2-2X_3X_1+2X_4-2t_2X_2X_1
+2(t_1+t_2)X_3 +2(t_1t_2+t_2^2)X_2; \end{eqnarray*}
and thus 
\begin{equation*}
R_1^q=R_1,\quad
R_2^{(q)} =  X_2^2-qX_1-2t_2X_2X_1
+(t_1+t_2)q
+2(t_1t_2+t_2^2)X_2.
 \end{equation*}

We define a graded $S[q]$-algebra by:
\[{\mathcal{B}}_n:=S[q][X_1,\cdots,X_{n}]/\langle 
R_1^{(q)},\ldots,R_n^{(q)} \rangle. \]
\begin{prop}\label{prop:A=B(typeC)}
We have an isomorphism of graded 
$S[q]$-algebras \[ \mathcal{B}_n\longrightarrow  \mathcal{A}_n\]
sending $X_i$ to $ Q_i(x|t)~(1 \le i \le n).$
\end{prop}
\noindent{Proof.} It suffices to show the analogue of Lemma \ref{lemma:univiso} .
Let 
\[
\widetilde{\mathcal{A}}_n:=
S[Q_1(x|t),\ldots,Q_{n}(x|t),2Q_{n+1}(x|t),Q_{n+2}(x|t),\ldots,Q_{2n}(x|t)].
\]
Consider the homomorphism of $S$-algebras
$$
\phi: 
\widetilde{\mathcal{B}}_n:=S[X_1,\ldots,X_n,2X_{n+1},X_{n+2},\ldots,X_{2n}]/
\langle R_1,\ldots,R_n\rangle
\rightarrow \widetilde{\mathcal{A}}_n.
$$ sending $X_i$ to $Q_i(x|t)$. We claim it is an isomorphism.
Surjectivity is obvious. Let us prove that $\phi$ is injective.
For arbitrary element $F$ in $\widetilde{\mathcal{B}}_n$, there is a sufficiently large integer $N$ such that 
$2^NF$ can be  represented by an element $G$ in $S[X_1,X_3,\ldots,X_{2n-1}]$ (cf. proof of Lemma \ref{lemma:univiso}).
Suppose $\phi(F)=0.$ Take $N$ and $G$ as above. Then $\phi(2^NF)=\phi(G)=0.$ Because $Q_{2i-1}(x|t)\;(1\leq i\leq n)$ are
algebraically independent over $S$, we have $G=2^NF=0$
in $\widetilde{\mathcal{B}}_n.$ 
A Gr\"obner basis argument similar to Lemma \ref{lem:free} implies that $\widetilde{\mathcal{B}}_n$  
is free over $S$ with an $S$-basis
$X^\la \cdot (2X_{n+1})^{m_{n+1}}\cdot X_{n+2}^{m_{n+2}}\cdots X_{2n}^{m_{2n}}$, where  $X^\lambda = X_{\lambda_1} \cdots X_{\lambda_k}$ for $\lambda= (\lambda_1, \ldots, \lambda_k) \in \mathcal{SP}(n)$ and $m_i\geq 0$. In particular, $\widetilde{\mathcal{B}}_n$ has no torsion element.
Hence we have $F=0$.
\qed

\bigskip

By the same argument of Lemma \ref{lem:free}, we see that $\mathcal{B}_n$ is free over $S[q]$
with basis $\{ X^\lambda\,|\;  \lambda \in \mathcal{SP}(n) \}$. We define the Pfaffian formula for $X_\la\in S[q][X_1,\ldots,X_n]$ analogous to (\ref{eq:Pfaff-P}): \begin{equation}\label{E:Xla}
X_\la=\mathrm{Pf} (X_{\la_i,\la_j})_{1\leq i<j\leq r}; \quad X_{k,l}=X_kX_l+\sum_{(r,s) \in \mathcal{I}_{k,l} }f_{k,l}^{r,s}(t)
X_rX_s 
\end{equation} with the substitutions (\ref{E:LGsub}) enforced. (We use again that the terms containing $X_{n+1}$ are divisible by $2$.) By the Pfaffian formula for the factorial $Q$-Schur functions it follows that $X_\lambda$ is sent to the image of $Q_\lambda(x|t)$ in $\mathcal{A}_n$
by the isomorphism of Proposition \ref{prop:A=B(typeC)}. The analogue of Proposition \ref{T:Ybasis} holds in this context as well with the same proof:

\begin{prop}\label{T:Xbasis} The polynomials $X_\lambda$, when $\lambda$ varies in $\mathcal{SP}(n)$, form a $S[q]$-basis for 
$\mathcal{B}_n.$ 
\end{prop}
We obtain the main result for the Lagrangian Grassmannian:

\begin{thm}\label{thm:LG} (i) (the factorial $Q$-Schur version) 
There is an isomorphism of 
graded $S[q]$-algebras \[
S[q][Q_1(x|t),\ldots,Q_n(x|t),2Q_{n+1}(x|t),Q_{n+2}(x|t),\ldots,Q_{2n}(x|t)]/J_n^{(q)}
\longrightarrow \QH_T^*(\LG(n)),
\]
where 
$J_n^{(q)}$ is the ideal generated by 
$2Q_{n+1}(x|t)-q,Q_{n+2}(x|t), \ldots,Q_{2n}(x|t)$.
Moreover, the image of $Q_\la(x|t)\;(\la \in \mathcal{SP}(n))$
is the Schubert class $\sigma_\la.$

(ii) (the indeterminate version) There is an isomorphism of 
graded $S[q]$-algebras
$$
S[q][X_1,X_2,\ldots,X_n]/\langle 
R_1^{(q)},R_2^{(q)},\ldots,R_n^{(q)}\rangle
\longrightarrow \QH_T^*(\mathrm{LG}(n))
$$
which sends $X_\la\;(\la \in \mathcal{SP}(n))$
to the Schubert class $\sigma_\la$. 
\end{thm}
\noindent{Proof.} The proof is the same as that of Theorem \ref{thm:OG}, using now the isomorphism $\mathcal{B}_n \simeq \mathcal{A}_n$ (Proposition \ref{prop:A=B(typeC)}), Proposition \ref{T:Xbasis}, and the Chevalley formula proved in Theorem \ref{T:EQCh-C}. \qed
 
%\bigskip
%\noindent{Remark 4.} The factorial $P$- and $Q$-functions satisfy certain left divided difference equations coming from the the left Weyl group action on the equivariant cohomology of either $\OG(n)$ or $\LG(n)$; see \cite{DSP}. Theorems \ref{thm:OG} and \ref{thm:LG} imply that  the equivariant quantum Schubert classes satisfy the same equations, and it is natural to ask whether these phenomena will hold for more general homogeneous spaces $G/P$. This issue will be discussed elsewhere.

%\setcounter{section}{1}
\setcounter{equation}{0}
%\addcontentsline{toc}{sectionlevel}{A}
\appendix
\renewcommand{\theequation}{A-\arabic{equation}}
\section{Equivariant cohomology of $\OG(n)$ via Chern classes}\label{ss:chernOG} %This section is  supplementary.
{By specializing Theorem \ref{thm:OG} at $q=0$,
we have a presentation for the equivariant cohomology ring of $\OG(n).$ 
%The goal of this section is to interpret the algebraic presentation for the equivariant cohomology of $\OG(n)$ in terms of the Chern classes of vector bundles. 
The goal of this section is to give a ``dictionary'' between the various algebraic
quantities in the presentation and the geometric quantities 
given in terms of equivariant Chern classes.}
In particular, we will show how the quadratic relations arise naturally from a Chern class calculation.
As we noted in the introduction, the formulas we obtain can be deduced from those for orthogonal degeneracy loci obtained of Kazarian \cite{kazarian}, although here we provide direct argument, similar to the one obtained by the authors in \cite[\S 11]{DSP} for $\LG(n)$.~In fact, this section can be seen as completing the aforementioned discussion from \cite{DSP}. 

\subsection{Presentation for $H_T^*(\OG(n))$}
{We start by recalling the following presentation for $H_T^*(\OG(n))$
obtained from the main theorem after specialization $q=0.$}
\begin{cor}\label{cor:presenHT} There is an isomorphism of 
graded $S$-algebras
\[
S[P_1(x|t),P_2(x|t),\ldots,P_{2n}(x|t)]
/I_n^0
\longrightarrow H_T^*(\OG(n)),
\]
where $I_n^0$ denote 
the ideal generated by 
$P_{n+1}(x|t),\ldots,P_{2n}(x|t).$
The map sends $P_\la(x|t)\;(\la\in \mathcal{SP}(n))$
to the equivariant Schubert class $\sigma_\la.$
\end{cor}

We leave the reader to write down the indeterminate version of this presentation. In that case, the ideal of relations
%{Clearly we also have the indeterminate version.
%We left the reader to write it down explicitly. Note that the ideal 
is generated by the elements obtained from the polynomials $R_i$ defined in (\ref{eq:R_i}) after substituting $X_i=0$ for $n+1\leq i\leq 2n$.
We will show in \S \ref{sec:quad} below why the quadratic relations $R_i=0$ are geometrically natural.

\bigskip 
\noindent{Remark 4.} Corollary \ref{cor:presenHT} can also be proved  
%We can prove Corollary \ref{cor:presenHT}
by using localization techniques.
%without using (quantum) equivariant cohomology.
In fact, we know that the factorial $P$-functions for $\la\in \mathcal{SP}(n)$ are sent to 
the equivariant Schubert classes in $H_T^*(\mathrm{OG}(n))$ % (\cite{EYD}, \cite{DSP})
by the map $\pi_n$ below (see (\ref{eq:piSchbert}) and Proposition \ref{prop:factorPi}).
Then an argument similar to the one in \cite[\S 8]{Ik} can be applied.

\bigskip

The Pfaffian formula (\ref{eq:Pfaff-P}) for $P_\la(x|t)$ implies
the following Giambelli formula
\begin{equation}
\sigma_\la=\mathrm{Pf}(\sigma_{\la_i,\la_j})_{1\leq i<j\leq r}.
\label{eq:eqGiam}
\end{equation}
This formula was proved in \cite{EYD} and 
also by Kazarian \cite{kazarian} in the context of degeneracy loci 
of vector bundles.

\bigskip
\subsection{Equivariant Chern classes, the map $\pi_n$}
%We will express the equivariant Schubert classes in terms of 
%the equivariant Chern classes of naturally defined
%vector bundles. 
 
We recall the setup from \cite[\S 10]{DSP}. There is %The form $\langle \cdot , \cdot \rangle_D$ determines 
a {\em tautological exact sequence\/} of bundles over $\OG(n)$ \begin{equation} 0 \to \mathcal{V} \to \mathcal{E} \to \mathcal{V}^* \to 0 \label{eq:sequence}\end{equation} where $\mathcal{E}$ is the trivial (but not equivariantly trivial) vector bundle $\C^{2n+2}$ and $\mathcal{V}$, the tautological subbundle,  has rank $n+1$. 
We identify $\mathcal{E}/\mathcal{V}$ with $\mathcal{V}^*$ by using the form $\langle \cdot , \cdot \rangle_D.$
The action of $T$ on $\C^{2n+2}$ determines a decomposition into weight spaces, which in turn determines a splitting $\mathcal{E} = (\oplus_{i=1}^{n+1} \mathcal{L}_i) \bigoplus (\oplus_{i=1}^{n+1} \mathcal{L}_i^*)$. Here $\mathcal{L}_i= \C {\pmb e}_i$, respectively $\mathcal{L}_i^* = \C {\pmb e}_i^*$ denote the trivial line bundle with $T$-weight $-t_i$, respectively $t_i$. Set $\mathcal{L} = \oplus_{i=1}^{n+1} \mathcal{L}_i$, $\mathcal{Z}_k = \oplus_{i=1}^k \mathcal{L}_i$ for $k \le n+1$,
and $ \mathcal{Z}_k =\mathcal{Z}_{n+1}$ for $k>n+1$.
If $E,F$ are vector bundles we denote by $c^T(E)$ the equivariant total Chern class of $E$ and by $c_i^T(E- F)$ the term of degree $i$ in the formal expansion of ${c^T(E)}/{c^T(F)}$. Define  \[ c_i^{(k)} =\textstyle{\frac{1}{2}}(c_i^T(\mathcal{V}^* - \mathcal{L} + \mathcal{Z}_k) - c_i^T(\mathcal{Z}_k)) \/. \] 
We have proved in \cite[Proposition 10.5]{DSP} that the class \[\gamma_i = \textstyle{\frac{1}{2}} c_i^T(\mathcal{V}^* - \mathcal{L})\] is an integral class in $H^*_T(\OG(n))$. The fact implies that $c_i^{(k)}$ is also integral. 

{
The sequence (\ref{eq:sequence}) leads to quadratic relations among $\gamma_i$'s 
$$
\gamma_i^2+2\sum_{j=1}^{i-1}(-1)^j\gamma_{i+j}\gamma_{i-j}+(-1)^{i}\gamma_{2i}=0 \quad (i\geq 1).
$$
One sees that the same relations are satisfied by 
$P_i(x)$'s (cf. Corollary \ref{cor:quadraticvanishing}). 
%\subsection{The map $\pi_n$}\label{sec:pi}
From this fact we can define} %In \cite[\S 10]{DSP} we constructed 
a $\Z[t]$-algebra homomorphism (\cite[Prop. 10.6]{DSP})
\[ \pi_n: \Z[t] \otimes_\Z \Gamma' 
\longrightarrow H^*_T(\OG(n)) \] such that $\pi_n(t_i) = 0 $ if $i > n+1$ and $\pi_n(P_i(x)) = \gamma_i$ for $i \ge 1 $.
The {fundamental} property satisfied by $\pi_n$ is that 
\begin{equation}
\pi_n(P_\lambda(x|t)) =  \begin{cases}
\sigma_\lambda & \lambda\in \mathcal{SP}(n)\\
0 & \lambda\not\in \mathcal{SP}(n).
\end{cases}\label{eq:piSchbert}\end{equation}
Consider the function $P_i^{(k)}(x|t)$ defined by (\ref{eq:2P_i}) above.
\begin{prop}\label{prop:pi_n(P)} We have
\begin{eqnarray}
\pi_n(P_i^{(k)} (x|t)) &=& \textstyle{\frac{1}{2}}c_i^T(\mathcal{V}^* - \mathcal{L} + \mathcal{Z}_k)\quad
(i\geq 1),\label{eq:piP1}\\
\pi_n(P_i(x|t))&=&c_i^{(i)} \quad(i\geq 1).\label{piP2}
\end{eqnarray}
\end{prop}
\noindent{Proof.} Straightforward by using the definitions and identities from \S \ref{ss:ChernPk}. 
\qed

\bigskip

The identities (\ref{eq:piSchbert}) and  (\ref{piP2}) imply:

\begin{cor}\label{cor:sigma_i} The class
$c_i^{(i)}$ is equal to $\sigma_i$ for $1\leq i\leq n$
and is equal to $0$ for $i>n.$
\end{cor}

\noindent{Remark 5.} Corollary \ref{cor:sigma_i} can also be proved by direct geometric arguments.~The equality $c_i^{(i)}= \sigma_i$ follows from the formulas for some special degeneracy loci obtained by Fulton and Pragacz - see e.g. \cite[pag. 90]{fulton.pragacz}. To prove the vanishing 
$c_i^{(i)}=0\;(i>n)$ first notice that if $i\geq n+1$ then  
$c_i^{(i)}=\frac{1}{2}(c_i^T(\mathcal{V}^*)-c_i^T(\mathcal{Z}_{n+1})).$ The vanishing for $i>n+1$ is 
a consequence that the bundles involved 
have rank $n+1$. Let us consider $i=n+1$ case. By definition of $\mathcal{V}$, 
the fibers of $\mathcal{V}$ and $\mathcal{Z}_{n+1}^*$ are in the same family, i.e. both fibers are on the same $\mathrm{SO}_{2n+2}$-orbit, 
and therefore also $\mathcal{V}^*$ and $\mathcal{Z}_{n+1}$ are. 
Then a result by Edidin-Graham \cite{edidin.graham:quadric} implies that $c^T_{n+1}(\mathcal{V}^*)=c^T_{n+1}(\mathcal{Z}_{n+1}),$
and the vanishing follows. 

\medskip

%We note also that this vanishing
%result explains why the presentation of Corollary \ref{cor:presenHT}
%exists from a geometric point of view.
\begin{prop}\label{prop:factorPi} The canonical projection 
$\pi_n$ factors as follows:
\[\Z[t]\otimes \Gamma'\longrightarrow
S[P_1(x|t),\ldots,P_{2n}(x|t)]/I_n^0
\overset{\cong}{\longrightarrow} H_T^*(\OG(n)),\]
where the second 
map is the isomorphism of  Corollary \ref{cor:presenHT}.
\end{prop}
\noindent{Proof.} This is clear from (\ref{eq:piSchbert}), or alternatively from (\ref{piP2}) and Corollary \ref{cor:sigma_i}. \qed
% implies that 
%$\pi_n(P_i(x|t))=0$ for $i>n.$ Then the result follows. 
%\qed

\bigskip

By applying $\pi_n$ to the equation from
Lemma \ref{lem:2row}, and using (\ref{eq:piSchbert}), (\ref{eq:piP1}), and (\ref{piP2}) one obtains:
\begin{prop}\label{prop:Chern2row}
For $1 \le l < k \le n$, we have \begin{equation} \sigma_{k,l} = c_k^{(k)} c_l^{(l)} + c_l^T(\mathcal{Z}_l) c_k^{(k)} + \sum_{j=1}^l (-1)^j c_{l-j}^T(\mathcal{V}^* - \mathcal{L} + \mathcal{Z}_l) c_{k+j}^{(k)}. \label{eq:tworowD} \end{equation}
\end{prop}

This finishes the interpretation in terms of the Chern classes of all the quantities involved in the presentation from Cororally \ref{cor:presenHT}.

\subsection{Quadratic relations}\label{sec:quad} We close this section by showing how the quadratic
relations for the factorial $P$-functions are naturally 
derived from geometric arguments.

\begin{prop}\label{prop:gquadratic} For each $1 \le k \le n$ we have the following relation in $H^*_T(\OG(n))$ \begin{equation}\label{E:gq} 
(c_k^{(k)})^2 + c_k^T(\mathcal{Z}_k) c_k^{(k)} + \sum_{i=1}^k (-1)^i c_{k-i}^T(\mathcal{V}^* - \mathcal{L} + \mathcal{Z}_k) c_{k+i}^{(k)} = 0 \/. \end{equation} \end{prop}

\noindent{Proof.} 
Let us denote the right hand side by $\mathcal{R}_k$. 
Using that $c^T(\mathcal{V} + \mathcal{V}^*) = c^T(\mathcal{E}) =  c^T(\mathcal{L} + \mathcal{L}^*)$ we obtain an equality 
\begin{equation} c^T(\mathcal{V}^* - \mathcal{L} + \mathcal{Z}_k) c^T(\mathcal{V} - \mathcal{L}^* + \mathcal{Z}_k^*)  = c^T(\mathcal{Z}_k + \mathcal{Z}_k^*) \/. \label{eq:quad-geom}
\end{equation} Notice that  $c_i^T(\mathcal{V} - \mathcal{L}^* + \mathcal{Z}_k^*) = (-1)^i c_i^T(\mathcal{V}^* - \mathcal{L} + \mathcal{Z}_k)$ and  $c_{k+i}^T(\mathcal{V}^* - \mathcal{L} + \mathcal{Z}_k) = 2 c_{k+i}^{(k)}$. Taking terms of degree $2k$ in both sides of (\ref{eq:quad-geom}), we obtain that $4 \mathcal{R}_k = 0$. Since the equivariant cohomology ring is torsion-free (being a free module over $S$), and $\mathcal{R}_k$ is an integral class, we have $\mathcal{R}_k = 0$. 
\qed

\begin{prop}[Geometric derivation of $R_k=0$] For $1\leq k\leq n,$
we have 
\begin{equation} P_k(x|t)^2 + (-1)^k e_k(t_1, \ldots, t_k) P_k(x|t) + 2 \sum_{i=1}^{k-1} (-1)^i P_{k-i}^{(k)}(x|t) P_{k+i}^{(k)}(x|t)  +  (-1)^k P_{2k}^{(k)}(x|t)=0. \label{eq:R=R'}
\end{equation}
\end{prop}
\noindent{Proof.}
Let us denote the left hand side by $R_k'.$
By Proposition \ref{prop:pi_n(P)}, it follows that  
$\pi_n(R_k')=\mathcal{R}_k.$
We know from \cite{DSP} that there is an injective $\Z[t]$-algebra homomorphism \[ \Z[t] \otimes_\Z \Gamma' \to \underleftarrow{\lim}_n H^*_T(\OG(n)).\] 
Note that the relations $\mathcal{R}_k =0$ also hold in $\underleftarrow{\lim}_n H^*_T(\OG(n))$ because they are compatible with the inverse system and they stabilize for large $n$. It follows that $R_k'=0$. 
Then it is straightforward to show  $R_k'=R_k$ by using Lemma \ref{lem:Pshit}.
\qed

\bigskip

\noindent{Remark 6.} The reader should have no trouble extending the arguments above for $\LG(n)$. Most of them are already present in \cite[\S 11.2]{DSP}. We only note that the analogue of the class $c_i^{(k)}$ in this case is $c_i^T(\mathcal{V}^* - \mathcal{L} + \mathcal{Z}_{k-1})$. 

\bigskip

Acknowledgements. We thank Anders Buch for helpful conversations related to this project.
T.I. was partially supported by Grant-in-Aid for
Scientific Research (C) 20540053; L.M. is partially supported by an NSA Young Investigator Grant H98230-13-1-0208;
H.N. is partially supported by Grant-in-Aid for
Scientific Research (C) 25400041.

\end{document}